\documentclass{article} 

\usepackage{amsmath,amssymb,amsthm}
\usepackage{graphicx,float,url,color,subcaption}
\usepackage{mathrsfs}
\usepackage[colorlinks=true]{hyperref}

\usepackage{calc}
\usepackage[usenames,dvipsnames]{xcolor}

\advance \topmargin by -\headheight
\advance \topmargin by -\headsep
     
\evensidemargin \oddsidemargin
\newcommand{\R}{\mathbb{R}}
\newcommand{\Z}{\mathcal{Z}}
\newcommand{\N}{\mathbb{N}}

\newcommand{\dx}{\Delta x}

\renewcommand{\P}{\mathcal{P}_c(\R^d)}
\newcommand{\Pac}{\mathcal{P}^{ac}_c(\R^{d})}
\renewcommand{\L}{\mathcal{L}}

\newcommand{\PR}{\mathcal{P}_c(\overline{B(0,R)})}
\newcommand{\PRac}{\mathcal{P}^{ac}_c(\overline{B(0,R)})}

\newtheorem{theorem}{Theorem}
\newtheorem{lemma}[theorem]{Lemma}
\newtheorem{corollary}{Corollary}
\newtheorem{prop}{Proposition}
\newtheorem{definition}{Definition}[section]
\theoremstyle{definition}\newtheorem{remark}{Remark}

\newcommand{\bt}{\begin{theorem}}
\newcommand{\et}{\end{theorem}}
\newcommand{\bl}{\begin{lemma}}
\newcommand{\el}{\end{lemma}}
\newcommand{\bp}{\begin{prop}}
\newcommand{\ep}{\end{prop}}
\newcommand{\bc}{\begin{corollary}}
\newcommand{\ec}{\end{corollary}}
\newcommand{\bdeff}{\begin{definition}}
\newcommand{\edeff}{\end{definition}}
\newcommand{\brem}{\begin{remark}}
\newcommand{\erem}{\end{remark}}

\newcommand{\eps}{\varepsilon}

\newcommand{\Pt}[1]{\left( #1 \right)}
\newcommand{\Pg}[1]{\left\{ #1 \right\}}
\newcommand{\Pq}[1]{\left[ #1 \right] }
\newcommand{\Pa}[1]{\langle #1 \rangle}
\newcommand{\Pabs}[1]{\left| #1 \right|}

\newcommand{\bi}{\begin{itemize}}
\newcommand{\ei}{\end{itemize}}

\newcommand{\be}{\begin{enumerate}}
\newcommand{\ee}{\end{enumerate}}

\newcommand{\bd}{\begin{description}}
\newcommand{\ed}{\end{description}}
\renewcommand{\i}{\item}

\newcommand{\bqn}{\begin{eqnarray}}
\newcommand{\eqn}{\end{eqnarray}}
\newcommand{\eqnn}{\nonumber\end{eqnarray}}
\newcommand{\eqnl}[1]{\label{#1}\end{eqnarray}}

\newcommand{\ba}[1]{\begin{array}{#1}}
\newcommand{\ea}{\end{array}}
\newcommand{\Lip}{\mathrm{Lip}}
\newcommand{\bproof}{\begin{proof}}
\newcommand{\eproof}{\end{proof}}

\renewcommand{\r}[1]{\eqref{#1}}

\newcommand{\supp}{\mathrm{supp}}
\newcommand{\weak}{\rightharpoonup}
\newcommand{\U}{\mathcal{U}}

\newcommand{\auth}[1]{{\sc #1}}
\newcommand{\tit}[1]{{\rm #1}}

\newcommand{\jou}[1]{{\it #1}}

\newcommand{\pp}[1]{pp.~#1}

\title{Mean-Field Sparse Jurdjevic--Quinn Control}

\author{Marco Caponigro\footnote{Conservatoire National des Arts et M\'etiers, \'Equipe M2N,
292 rue Saint-Martin, 75003, Paris, France ({\tt\small marco.caponigro@cnam.fr}).}, 
\and
Benedetto Piccoli\footnote{Department of Mathematical Sciences and Center for Computational and Integrative Biology, Rutgers University, Camden, NJ 08102, USA ({\tt\small piccoli@camden.rutgers.edu}).}
\and
Francesco Rossi\footnote{Aix Marseille Universit\'e, CNRS, ENSAM, Universit\'e de Toulon, LSIS, Marseille, France ({\tt\small francesco.rossi@lsis.org}.)
}
\and
Emmanuel Tr\'elat\footnote{Sorbonne Universit\'es, UPMC Univ. Paris 6, CNRS UMR 7598, Laboratoire Jacques-Louis Lions, 4 place Jussieu, 75005, Paris, France ({\tt\small emmanuel.trelat@upmc.fr}).}}

\begin{document}

\maketitle

\begin{abstract}
We consider nonlinear transport equations with non-local velocity, describing the time-evolution of a measure, which in practice may represent the density of a crowd. 
Such equations often appear by taking the mean-field limit of finite-dimensional systems modelling collective dynamics.
We first give a sense to dissipativity of these mean-field equations in terms of Lie derivatives of a Lyapunov function depending on the measure.

Then, we address the problem of controlling such equations by means of a time-varying bounded control action localized on a time-varying control subset with bounded Lebesgue measure (sparsity space constraint). Finite-dimensional versions are given by control-affine systems, which can be stabilized by the well known Jurdjevic--Quinn procedure.

In this paper, assuming that the uncontrolled dynamics are dissipative, we develop an approach in the spirit of the classical Jurdjevic--Quinn theorem, showing how to steer the system to an invariant sublevel of the Lyapunov function. The control function and the control domain are designed in terms of the Lie derivatives of the Lyapunov function, and enjoy sparsity properties in the sense that the control support is small.

Finally, we show that our result applies to a large class of kinetic equations modelling multi-agent dynamics. 
\end{abstract}


\section{Introduction and main result}
\subsection{The context}
In recent years, the study of collective behavior of a crowd of autonomous agents has drawn a great interest from scientific communities, e.g., in civil engineering (for evacuation problems \cite{CPT-14,HFV-00}), robotics (coordination of robots \cite{BCM,JLM-03,OSFM,RME}), computer science and sociology (social networks \cite{HG-02}), and biology (animals groups \cite{BB-11,CKJRF,G-08}). In particular, it is well known that some simple rules of interaction between agents can promote formation of special patterns, like lines in ants formations and migrating lobsters, or V-shaped formation in migrating birds. This phenomenon is often referred to as {\em self-organization}.

Beside the problem of analyzing the collective behavior of a ``closed'' system~\cite{CFP11}, it is interesting to understand what changes of behavior can be induced by an external agent (e.g., a policy maker) to the crowd. For example, one can try to enforce creation of patterns when they are not formed naturally, or break the formation of such patterns \cite{CS-control,CS-control-1,FPR-14,controlKCS, HPS15}. This is the problem of {\em control of crowds}, that we address in this article in a specific case.

From the mathematical point of view, problems related to models of crowds are of great interest. From the analysis point of view, one needs to pass from a big set of simple rules for each individual to a model capable of capturing the dynamics of the whole crowd. This can be solved via the so-called mean-field process, that permits to consider the limit of a set of ordinary differential equations (one for each agent) to a partial differential equation (PDE in the following) for the whole crowd
\cite{controlKCS}.
The resulting equation is a transport equation with non-local velocity, of the form
\bqn
\partial_t \mu + \nabla\cdot (f\Pq{\mu}\mu)=0,
\eqnl{e-pdebase}
where $\mu$ is the measure representing the density of agents, $\nabla\cdot$ is the divergence operator and $f[\mu]$ is a vector field depending on the measure, taking into account interactions between agents. Many kinetic equations are of this form.
We recall fundamental properties of such classes of equations in Section \ref{s-PDE}. We just highlight here that we assume in the following that $f[\mu]$ is a bounded Lipschitz vector field for any $\mu$, Lipschitz with respect to the Wasserstein distance $W_p$,  $p\in[1,+\infty)$,  as a function of $\mu$. This ensures existence and uniqueness of the solution of the associated Cauchy problem \cite{ambrosio-gangbo,pedestrian}. Since such an equation generates a semigroup, we will use the notation $e^{t f}\mu_0$ to denote the unique solution of \r{e-pdebase} at time $t$ with initial data $\mu_0$. We recall existence results for \r{e-pdebase} in Section~\ref{s-PDE}.

\subsection{Control of transport equations with non-local velocity}
To control systems described by \r{e-pdebase}, we assume to act only on a small part of the crowd. Since agents are indistinguishable when one only knows $\mu(t)$ at time $t$, controls can only be state-dependent and cannot focus on specific agents. For this reason, we model the control action by means of a vector field $g$, and a {\em control gain} $u(t,x)$ localized in a small {\em control set} $\omega(t)$ (itself depending on time), modeling our choice of the gain on the vector field. We choose the control gain $u$ and the control set $\omega$, both of them varying in time while the vector field $g$ is fixed, and it depends on the density $\mu$. The resulting control system is given by
\bqn
\partial_t \mu +\nabla\cdot ((f\Pq{\mu}+  \chi_{\omega} u g[\mu]) \mu)=0.
\eqnl{e-PDE}
Here, the function $\chi_\omega$ is the indicator function of $\omega$, defined almost verywhere by $\chi_\omega(x)=1$ if $x\in\omega$ and $\chi_\omega(x)=0$ otherwise.

We focus on the following {\em continuous sparse space constraint}: we assume to act only on a small portion of the configuration space, with finite strength. Accordingly, we assume the following constraints. Here, given a measurable subset $\omega$ of $\R^d$, we denote by $|\omega|$ its Lebesgue measure.
\\

\noindent\fbox{%
    \parbox{0.99\textwidth}{%
    \begin{center}
    {\bf Control Constraints (U)}
    \end{center}

Fix $c>0$. For each time $t\geq0$, we have:
\bqn
&\mbox{{\bf Sparsity space constraint:}} &|\omega(t)| \leq c,\label{e-spacecon}\\
&\mbox{{\bf Finite strength:}} & \| u(t,\cdot)\|_{L^\infty}\leq 1.\label{e-strength}
\eqn

}}

\vspace{3mm}

Control sparsity constraints have been first introduced in \cite{CS-control,CS-control-1}, for a population with a finite number of agents. The sparsity space constraint was considered in \cite{controlKCS}. In the mean-field approach, this is actually the most natural sparsity constraint when one wants to use space-dependent vector fields and to act on ``small sets''. A {\em sparsity population constraint} was also considered in \cite{controlKCS}.

For the sparsity space constraint, one can easily deal both with measures that are absolutely continuous with respect to the Lebesgue measure, and with measures containing singular (Dirac) parts. We denote by $\P$ the space of probability measures on $\R^d$ with compact support and by $\Pac$ the subspace of probability measures on $\R^d$ that are absolutely continuous with respect to the Lebesgue measure.

In the following, we will define a control strategy that satisfies the following property: {\em if the initial data, at time $0$, is absolutely continuous} with respect to the Lebesgue measure, {\em then it remains absolutely continuous for any positive time}. This does not prevent $\mu(t)$ of converging to some Dirac mass as $t\rightarrow+\infty$, as this is the case for consensus problems for crowd models.

\medskip

In this paper, our objective is to generalize the Jurdjevic--Quinn stabilization method \cite{JQ} to mean-field controlled equations, under the sparsity constraint {\bf (U)} described above. Following the Jurdjevic--Quinn approach, we assume to have available a Lyapunov function $V$ for which:
\begin{itemize}
\i the uncontrolled dynamics $f[\mu]$ gives no increase of $V$;
\i the control $u g[\mu]$ allows one to increase-decrease $V$, except for some specific configurations of the population $\mu$ in the subset $\mathcal{Z}$ of the set of measures with compact support $\P$ defined as the set on which the Lie derivatives of $V$ vanish (see the precise definition in~\eqref{e-Z}).
\end{itemize}

We will then define a sparse control strategy, steering the population exactly to the set $\mathcal{Z}$, in complete analogy with the standard finite-dimensional Jurdjevic--Quinn method.

The vector fields $f$ and $ug$ defined on the space $\R^d$ play the role of derivatives for the function $V[\mu]$, in the following sense: the vector field $f$ induces an infinitesimal change in $V$ that can be estimated as the derivative $\lim_{t\to 0}\frac{V[e^{t f}\mu]-V[\mu]}{t}$. This limit needs to be well-defined, and this is why we will assume from now on the following regularity assumptions on $f$: for any $\mu\in\P$, the function $t\mapsto V[e^{tf}\mu]$ is of class $C^2$. Then, in analogy with the definition of a Lie derivative in finite dimension, we give the following definition.

\begin{definition}
We define the \emph{Lie derivative} of $V$ along $f$ as the limit
\bqn
\L_{f} V[\mu]=\lim_{t\to 0}\frac{V[e^{t f}\mu]-V[\mu]}{t}.
\eqnl{e-Lie}
\end{definition}

Requiring the non-increase of $V$ along the flow of $f$ is equivalent to require dissipativity for the system.

\begin{definition}
We say that the system~\eqref{e-PDE} is dissipative if there exists a Lyapunov function $V:\P \to \R$ such that  $t\mapsto V[e^{tf}\mu]$ is of class $C^2$  and
\bqn
\L_{f}V[\mu]\leq 0\qquad\forall\mu\in\P.
\eqnl{e-dissipative}
\end{definition}

Regularity conditions need to be satisfied as well for the controlled vector field $u g$. For this reason, we define the space $\U$ of admissible control functions $u$, and from now on, we  assume that $\U=\Lip(\R^d,\R)$. In what follows, we impose the following regularity assumption: for all $\mu\in\P$ and $u\in \U $, the function $t\mapsto V[e^{tug}\mu]$ is of class $C^2$. As a consequence, the limit
\bqn
\L_{ug} V[\mu]=\lim_{t\to 0}\frac{V[e^{t ug}\mu]-V[\mu]}{t}
\eqnl{e-Lie2}
is well defined.

\begin{remark}
Accordingly, the notion of Lie derivative can be extended to piecewise constant functions 
$t \mapsto u(t,\cdot) \in \U$.
\end{remark}

\begin{definition}
Assume that $V:\P \to \R$ is such that  $t\mapsto V[e^{tf}\mu]$ and $t\mapsto V[e^{tug}\mu]$ are of class $C^2$ for all $\mu\in\P$ and $u\in \U$. We define
\bqn
\mathcal{Z}=\Pg{\mu\in\P\ \mid\ \L_fV[\mu]=\L_{ug}V[\mu]=0\quad\forall u\in\Lip(\R^d,\R)}.
\eqnl{e-Z}
\end{definition}

Note that the definition of Lie derivative implies the multiplicative property
\bqn
\L_{k f} V[\mu]=k\L_{f}V[\mu]\mbox{~~~~and~~~~}\L_{k ug} V[\mu]=k\L_{ug}V[\mu].
\eqnl{e-mult}
This continuity condition also implies additivity of Lie derivatives. Indeed, one can easily see that $e^{t(u+u')g}=e^{tug+o(t)}e^{tu'g}$, which in turn implies that
\bqn
\L_{(u+u')g} V[\mu]=\L_{ug}V[\mu]+\L_{u'g}V[\mu].
\eqnn
While conditions \r{e-Lie}-\r{e-Lie2} are equivalent to differentiability of $V$ along $f$ and $ug$, we also need a kind of differentiability for $V$ along directions of the dynamics. By the additivity property, we can state it as follows: there exists $K>0$ such that, for all $\mu\in\P$ and $u\in \U$, we have
\bqn
\Pabs{\L_{ug}V[\mu]}\leq K \|u\|_{L^1(\mu)}.
\eqnl{e-diff}
This yields a metric for the space of controls $u$ similar to the zero-order metric in a more general sub-Riemannian structure for metrics on the space of diffeomorphisms on a manifold \cite{arg-tre} (see also \cite{agra-capo}). The main difference here is that we choose the $L^1$ norm weighted with respect to the measure $\mu(t)$, and not with respect to the Lebesgue measure. 

While conditions \r{e-Lie}-\r{e-Lie2}-\r{e-diff} hold for a fixed $\mu\in\P$, we also require the continuity of the Lie derivatives of first and second order. We require that
\bqn
\ba{l}
\displaystyle\lim_{i\to+\infty}\L_{f}V[\mu^i]=\L_{f}V[\mu],\qquad \quad
\displaystyle\lim_{i\to+\infty}\L_f^2V[\mu^i]=\L_f^2V[\mu],\\
\displaystyle \lim_{i\to+\infty}\L_{ug}V[\mu^i]=\L_{ug}V[\mu],\qquad
\displaystyle \lim_{i\to+\infty}\L_{ug}^2V[\mu^i]=\L_{ug}^2V[\mu],\\
\displaystyle \lim_{i\to+\infty}\L_{f}\L_{ug}V[\mu^i]=\L_{f}\L_{ug}V[\mu] ,
\ea
\eqnl{e-contf}
for all $\mu^i$ and $\mu\in\P$, with $\mu^i\weak \mu$ (weak convergence of measures), i.e., $\lim_{i\to+\infty}\int \phi \,d\mu^i=\int \phi\,d\mu$ for every $\phi\in C^0_c(\R^d,\R)$. The conditions \r{e-contf} imply in particular that $\L_{f+ug}V[\mu]$ exists and satisfies $\L_{f+ug}V[\mu]=\L_fV[\mu]+\L_{ug}V[\mu]$, i.e., that additivity holds also for the vector field $f+ug$.

\brem Clearly, the choice of the set of admissible controls $\U$ has an impact on the set of admissible functionals $V$ for which \r{e-Lie2}-\r{e-diff}-\r{e-contf} are satisfied. We choose here the set of Lipschitz functions because existence is then ensured for \r{e-PDE} (see \cite{ambrosio-gangbo,pedestrian}).

Note that reducing the space of admissible controls to some proper subset of $\Lip(\R^d,\R)$ (such as $C^\infty_c(\R^d,\R)$) may enlarge the set of functionals $V$ for which the regularity conditions \r{e-Lie}-\r{e-Lie2}-\r{e-diff}-\r{e-contf} are satisfied. In Section~\ref{s-concentration}, we will enforce the decrease of the functional $V$ by a steepest descent method on the space $\U$, by (approximately) solving an optimization problem in the space of Lipschitz functions.
\erem 

Following the classical Lyapunov theory for finite-dimensional systems, we need to impose some conditions ensuring compactness of trajectories. In finite dimension, this is often imposed by requiring $V$ to be proper, i.e., $\lim_{\vert x\vert\to+\infty}V(x)=+\infty$, hence the fact that $\frac{d}{dt} V(x(t))\leq 0$ implies compactness. In the present mean-field setting, instead, such a condition cannot be imposed by a simple evaluation of the function $V$: in the case where $V$ is the variance of the measure, one can have measures $\mu$ with arbitrarily small variance and arbitrarily large support.
For this reason, we impose compactness of trajectories by assuming that the dynamics of the system, i.e., the vector fields $f$ and $g$, have a compact support. More precisely, we assume the existence of a ball $B(0,R)$ such that 
\begin{equation}
\mbox{ for all } \mu \in \PR, \qquad  f[\mu](x)=g[\mu](x)=0\qquad\forall x\not\in\overline{B(0,R)}.\label{e-compact}
\end{equation}
Note that this condition implies that $\mu(t)\in\PR$ for any $t\geq 0$. Then, all assumptions described above will be assumed to hold for measures in $\PR$ only.

\begin{remark}
 The results of this paper can be stated for transport equations with non-local terms~\eqref{e-PDE} defined on bounded manifolds without boundary, or on bounded manifold with no-flux boundary condition.
\end{remark}

Summing up, we make the following assumptions on the system \r{e-PDE} and on the Lyapunov functional $V$:\\

\noindent\fbox{%
    \parbox{0.99\textwidth}{%
    \begin{center}
    {\bf Assumptions (H)}
    \end{center}
    
The vector fields $f,g:\P\to \mathrm{Lip}(\R^d,\R^d)$ and the Lyapunov function $V:\P \to \R$ satisfy Assumptions {\bf (H)} if:
    \bi
     \i there exists $R>0$ such that $f,g$ satisfy the compact support property \r{e-compact};
    \i there exist $L>0$, $Q>0$ and $p\geq 1$ such that, for all $\mu,\nu\in\PR$ and for all $x,y\in\R^d$, 
\begin{equation}\label{e-regPDE}
\begin{split}
& |f[\mu](x)-f[\mu](y)|\leq L |x-y|,\qquad |g[\mu](x)-g[\mu](y)|\leq L |x-y|, \\
& |f[\mu](x)-f[\nu](x)|\leq Q W_p(\mu,\nu), \quad |g[\mu](x)-g[\nu](x)|\leq Q W_p(\mu,\nu) ;
\end{split}
\end{equation}
 \i the functions $t\mapsto V[e^{tf}\mu]$ and $t\mapsto V[e^{tug}\mu]$ are of class $C^1$ for all $\mu\in\PR$ and $u\in\U$; in particular, the Lie derivatives \r{e-Lie} and \r{e-Lie2} exist;
 \i the uncontrolled system is dissipative, i.e., $\L_{f}V[\mu]\leq 0$ for any $\mu\in \PR$;
 \i the Lie derivative \r{e-Lie2} satisfies the Lipschitz condition \r{e-diff} and the continuity condition \r{e-contf} for any $\mu\in \PR$.
    \ei
 }}

\begin{remark}\label{rk:bound}
The uniform Lipschitz property of $f$ and $g$ in~\eqref{e-regPDE} and the uniform compactness of their support in~\eqref{e-compact} imply that there exists $M>0$ such that $\|f[\mu]\|_{L^\infty}\leq M$ and $\|g[\mu]\|_{L^\infty}\leq M$ for any $\mu\in \mathcal{P}_c(\overline{B(0,R)})$.  
These facts imply existence and uniqueness of the solution of the Cauchy problem for \r{e-PDE} (see, e.g., \cite{ambrosio-gangbo,pedestrian}). We recall existence results in Section \ref{s-PDE}. 
\end{remark}

\subsection{The main result}
The main idea of our control strategy is to choose the controller to make $V$ decrease along trajectories. We will do this choice with a steepest descent method, similarly to the finite-dimensional approach described in \cite{CS-control,CS-control-1,finiteJQ}.

Since the space of admissible controls $\chi_{\omega} u g[\mu]$ is infinite dimensional, we restrict ourselves to a finite-dimensional set by imposing the following structure.
Consider the class of Lipschitz mollified indicator functions $\chi_{[a,b]}^\eta:\R\to\R$, defined by
\bqn
\chi^\eta_{[a,b]}(x)=
\begin{cases}
1&\mbox{ for }x\in[a,b],\\
0&\mbox{~for } x\not\in[a-\eta,b+\eta],\\
\frac{x-a+\eta}{\eta}& \mbox{~for } x\in[a-\eta,a];\\
\frac{-x+b+\eta}{\eta}& \mbox{~for } x\in[b,b+\eta],
\end{cases}
\eqnn
and then, consider the $d$-dimensional version of such functions. Given $a = (a^1,\ldots,a^d), b= (b^1,\ldots,b^d)$ and $x = (x^1,\ldots,x^d) $ in $\R^d$, we define
\bqn
\chi^\eta_{[a,b]}(x)= \min_{i=1,\ldots,d}  \chi^\eta_{[a^i,b^i]}(x^i) .
\eqnl{e-chi}
Now, for any choice of the three parameters $(a,b,\eta)$, we take $\omega=\omega(a,b,\eta)$ as  the multi-interval $[a^1-\eta,b^1+\eta]\times [a^2-\eta,b^2+\eta]\times\ldots\times [a^d-\eta,b^d+\eta]$. 
Then we reduce the choice of the sparse control in an infinite-dimensional space of controls 
 to the choice of three parameters $(a,b,\eta)$. 
In what follows, we set 
\begin{equation}\label{e-U}
 U(a,b,\eta) 
 =\chi_{[a^1-\eta,b^1+\eta]\times [a^2-\eta,b^2+\eta]\times\cdots\times [a^d-\eta,b^d+\eta]}  \chi^\eta_{[a,b]}.
\end{equation}
We then define the ``slope function'' by
$$
s_{t}(a,b,\eta)=|\L_{U(a,b,\eta)g[\mu(t)]} V[\mu(t)]|,
$$
which describes the instantaneous variation of $V$ in $\mu(t)$ as a consequence of the action of the control $U(a,b,\eta)$. Note that \r{e-diff} and the fact that the function $t\mapsto V[e^{tug}\mu]$ is of class $C^1$ imply the continuity of the slope function with respect to its arguments $(t,a,b,\eta)$.

We then apply a steepest descent method by choosing the control corresponding to one of the maximizers\footnote{The method to select a maximizer plays no role in the convergence of the method. One may consider for instance the lexicographic order in $\R^{2n+1}$, and choose the smallest maximizer.} $(a^*,b^*,\eta^*)$ of $s_t$  in the space 
\bqn
\Omega_t=\Pg{(a,b,\eta)\ \mid\ |\omega(a,b,\eta)|\leq c \mbox{  and } \eta\geq t^{-1}}.
\eqnl{e-Omega}
The condition $|\omega(a,b,\eta)|\leq c$ in \r{e-Omega} ensures that the space constraint \r{e-spacecon} is satisfied. We will see in Lemma~\ref{l-Linf} that the condition $\eta\geq t^{-1}$ implies that the control function is uniformly Lipschitz for any bounded time interval $[0,\theta]$, thus ensuring that $\mu(t)$ remains absolutely continuous with respect to the Lebesgue measure. At the same time, when $t\to+\infty$, this constraint allows to consider controls with an arbitrarily large Lipschitz constant, since $\Lip(\chi^\eta_{[a,b]})=\frac1\eta$. This Lipschitz constraint is somehow unavoidable if one wants to ensure regularity of the measure $\mu(t)$ within finite time; otherwise, the steepest descent method might either generate a non-Lipschitz vector field (for which existence for \r{e-PDE} holds for small times only) or a time-varying Lipschitz vector field converging to a non-Lipschitz vector field within finite time (see an example for a problem of crowd dynamics in Section~\ref{s-concentration}).

Choosing the control as the instantaneous maximizer of $s_t$ may cause chattering (in time) phenomena, as it has been already noticed in finite dimension (see \cite{finiteJQ}). For this reason, we regularize the control by means of an {\it hysteresis}: we introduce a parameter $h\in(0,1)$ and, given the control $U(a^*,b^*,\eta^*)$, maximizer of $s_t$ at time $t_n$, we keep it constant over an interval $[t_n,t_n+\delta]$ along which 
$(s_{t}(a^*,b^*,\eta^*)\geq (1-h)s_{t}(a,b,\eta)$. 

Summing up, the combination of a {\em steepest descent method} with an {\em hysteresis} provides a control making $V$ decrease and steering the density $\mu(t)$ to $\mathcal{Z}$. Our main result is the following.

\bt[Main Theorem] \label{t-main}
Let $f, g: \P\to \mathrm{Lip}(\R^d,\R^d)$ and $V:\P\to \R$ satisfy Assumptions {\bf (H)} for some $R>0$.
Consider the controlled transport equation with non-local velocity
\bqn
\partial_t \mu +\nabla\cdot \Pt{(f\Pq{\mu}+\chi_{\omega(t,y)} u(t,y) g[\mu]) \mu}=0,\qquad
\mu(0)=\mu_0,
\eqnl{e-PDEt}
where $\mu_0\in\Pac$ is such that $\supp(\mu_0)\subset \overline{B(0,R)}$.
Fix the hysteresis parameter $h\in (0,1)$. Fix the following initial parameters $n=0$ and $t_0=0$. 
Define the following algorithm step.\\

\noindent\fbox{%
    \parbox{0.99\textwidth}{%
    
        \begin{center}
    {\bf Step $n$}
    \end{center}

At time $t_n$, choose one of the maximizers $(a^*,b^*,\eta^*)$ of $s_{t_n}(a,b,\eta)$ in the set $\Omega_{t_n}$ defined in \r{e-Omega}.

Then, we have two cases:
\bi
\i If either $s_{t_n}(a^*,b^*,\eta^*)< t_n^{-1}$ or $\Omega_{t_n}$ is empty, then choose the zero control 
\bqn
\chi_{\omega(t)} u(t,x)\equiv 0 ,
\eqnl{e-control0} 
(thus, $\omega(t)$ needs not be defined)
and let the measure $\mu(t)$, starting at $\mu(t_n)$, evolve according to \r{e-PDEt} over the time interval $[t_n,t_{n+1}]$, where $t_{n+1}$ is the smallest time greater than $t_n$ for which
there exists $(\bar a,\bar b,\bar \eta)\in\Omega'_{t}$ such that $s_{t_{n+1}} (\bar a,\bar b,\bar \eta) \geq 2t_{n+1}^{-1}$, where 
\bqn
\Omega'_t=\Pg{(a,b,\eta)\in\Omega_t\ \mid\ \eta\geq 2t^{-1}}.
\eqnl{e-Omega1}

\i If $s_{t_n}(a^*,b^*,\eta^*)\geq t_n^{-1}$, then choose the control defined by
\bqn
\begin{split}
\omega(t) &= [{a^*}^1-\eta^*,{b^*}^1+\eta^*]\times [{a^*}^2-\eta^*,{b^*}^2+\eta^*]\times\cdots\times [{a^*}^d-\eta^*,{b^*}^d+\eta^*] , \\
u(t,\cdot)&= -\chi_{[a^*,b^*]}^{\eta^*}
\mathrm{sign}(\L_{U(a,b,\eta)g[\mu(t)]} V[\mu(t)]) ,
\end{split}
\eqnl{e-control}
where $U$ is given in ~\eqref{e-U}, and let the measure $\mu(t)$, starting at  $\mu(t_n)$, evolve according to \r{e-PDEt} over the time interval $[t_n,t_{n+1}]$, where $t_{n+1}$ is the smallest time greater than $t_n$ satisfying at least one of the following conditions:
\bi
\i either $s_{t_{n+1}}(a^*,b^*,\eta^*)\leq \frac{t_{n+1}^{-1}}2$;
\i or there exists $(\bar a,\bar b,\bar \eta)\in \Omega'_{t_{n+1}}$ such that 
\bqn
s_{t_{n+1}}(a^*,b^*,\eta^*)\leq (1-h)s_{t_{n+1}} (\bar a,\bar b,\bar \eta). 
\eqnl{e-hyst}
\ei
\ei
If $t_{n+1}$ is finite, then go to Step $(n+1)$. 

If $t_{n+1}=+\infty$, then keep the control \eqref{e-control0} or \eqref{e-control} over the time interval $[t_n,+\infty)$.
}}

\vspace{3mm}

For this control strategy, the control $\chi_{\omega} u$ satisfies the control constraint {\bf (U)}, the unique solution $\mu(t)$ of \eqref{e-PDEt} is such that $\mu(t)\in \Pac$ for any $t\in[0,+\infty)$, and $\mu(t)$ converges
to $\mathcal{Z}\cap\PR$, i.e. 
\bi
\i $\lim_{t\to+\infty}\inf_{\nu\in\mathcal{Z}\cap\PR}W_p(\mu(t),\nu)=0,$ 
\i or equivalently, there exists a choice $\nu(t)\in\mathcal{Z}\cap\PR$ for each $t\geq 0$ such that $\mu(t)\weak\nu(t)$, i.e. for all $\phi\in C^\infty_c(\R^d)$ it holds $\lim_{t\to\infty} \int\phi d(\mu(t)-\nu(t))=0$.
\ei
\et

\begin{remark} The three threshold time-dependent functions used in the definition of control algorithm in Theorem~\ref{t-main} satisfy $t^{-1}/2<t^{-1}<2t^{-1}$. One can easily see that they can be replaced with three positive functions satisfying $\phi_1(t)<\phi_2(t)<\phi_3(t)$ converging to $0$ as $t\to+\infty$. In particular, the functions can take finite values for $t=0$, by maybe allowing one control to be active on the starting interval $[0,t_1]$.
\end{remark}

The rest of the paper is structured as follows. In Section \ref{s-PDE}, we recall the main definitions and results for transport PDEs with non-local velocities as \r{e-pdebase} and \r{e-PDE}. In Section \ref{s-concentration}, we discuss some examples of dynamics of the form \r{e-PDE}, and we show explain some differences with respect to the finite dimensional setting. Theorem \ref{t-main} is proved in Section \ref{s-proof}. In Section~\ref{sec:gen}, we study a generalization of the Theorem~\ref{t-main} to a system of the form~\eqref{e-PDE}  with several control potentials. Finally, in Section~\ref{sec:mas}, we present an application of Theorem~\ref{t-main} to the control of kinetic multi-agent systems.

\section{Transport equations with non-local velocities}
\label{s-PDE}

In this section, we recall existence and uniqueness results for \eqref{e-pdebase} and \r{e-PDE}.
In \r{e-pdebase}, the variable $\mu\in \P$ is a probability measure on $\R^d$. The term $f[\mu]$ is called the \textit{velocity field} and it is a non-local term. Since the value of a measure at a single point is not well defined, it is important to observe that $f[\mu]$ is not a function depending on the value of $\mu$ in a given point, as it is often the case in the setting of hyperbolic equations in which $f[\mu](x)=f(\mu(x))$. Instead, one has to consider $f$ as an operator taking an as input the whole measure $\mu$ and giving as an output a global vector field $f[\mu]$ on the whole space $\R^d$. These operators are often called ``non-local'', as they consider the density not only at a given point, but in a whole neighbourhood.

We first recall two useful definitions to deal with measures and solutions of \eqref{e-pdebase}, namely the Wasserstein distance and the push-forward of measures (for more details see, e.g., \cite{villani}).

\begin{definition}
Given two probability measures $\mu$ and $\nu$ on $\R^d$ and $p\in[1,+\infty)$, the $p$-Wasserstein distance between $\mu$ and $\nu$ 
is
$$
W_p(\mu,\nu) = \inf\left\{ \int_{\R^{2d}} |x-y|^p\,d\pi(x,y) \mid  \pi\in\Pi(\mu,\nu) \right\}^{1/p},
$$
where $\Pi(\mu,\nu)$ is the set of transference plans from $\mu$ to $\nu$, i.e., of the probability measures $\pi$ on $\R^{d} \times \R^d$ such that 
$\mathrm{Proj}_x \# \pi =\mu$ and $\mathrm{Proj}_y \#\pi=\nu$ with  $\mathrm{Proj}_x:(x,y) \mapsto x$ and $\mathrm{Proj}_y:(x,y) \mapsto y$.
\end{definition}

The topology induced by $W_p$ on the space of probability measures $\mathcal{P}(X)$ on a compact space $X$ coincides with the weak-$*$ topology of measures (see \cite[Theorem 7.12]{villani}). As a consequence of  condition \r{e-compact}, each trajectory $\mu(t)$ of the controlled system \r{e-PDE} is contained in the compact space $\PR$ (compact if endowed with the Wasserstein topology). Thus, {\bf from now on, we will state equivalently convergence with respect to the  weak-$*$ topology of measures and with respect to the Wasserstein distance}.
We now define the push-forward of measures.

\begin{definition} 
Given a Borel map $\gamma:\R^d\rightarrow\R^d$, the pushforward of a measure $\mu\in \P$ is defined by
$\gamma\#\mu(A)=\mu(\gamma^{-1}(A))$
for every measurable subset $A$ of $\R^d$.
\end{definition}

We now recall an existence and uniqueness result for \eqref{e-pdebase} (see a complete proof in  \cite{controlKCS}).

\begin{theorem}\label{t-esistenza}
We assume that, for every $\mu\in \mathcal{P}_c(\R^d)$, the velocity field $f[\mu]$ is a function of $(t,x)$ with the regularity
\begin{equation*}
\begin{array}{rcl}
f[\cdot]: \mathcal{P}(\R^d) & \longrightarrow & \Lip(\R^d)\cap L^\infty(\R^d) \\
\mu & \longmapsto & f[\mu] 
\end{array}
\end{equation*}
satisfying the following assumptions:
\begin{itemize}
\item 
there exist functions $L(\cdot)$ and $M(\cdot)$ in $L^\infty_{loc}(\R)$ such that 
$$
\left\Vert f[\mu](t,x)-f[\mu](t,y) \right\Vert \leq L(t) \Vert x-y\Vert,\qquad \Vert f[\mu](t,x)\Vert \leq M(t) (1+\Vert x\Vert ),
$$
for every $\mu\in\mathcal{P}_c(\R^d)$, every $t\in\R$ and all $(x,y)\in\R^d\times\R^d$;
\item for a given $p\in[1,+\infty)$, there exists a function $K(\cdot)$ in $L^\infty_{loc}(\R)$ such that  
$$
\left\Vert f[\mu]-f[\nu] \right\Vert _{L^\infty(\R;C^0(\R^d))} \leq K(t) W_p(\mu,\nu),
$$
for all $(\mu,\nu)\in(\mathcal{P}_c(\R^d))^2$.
\end{itemize}
Then, for every $\mu^0\in\mathcal{P}_c(\R^d)$, the Cauchy problem
\begin{equation}\label{e-cauchy}
\partial_t\mu+\nabla\cdot(f[\mu] \mu)=0, \quad \mu_{|_{t=0}}=\mu_0,
\end{equation}
has a unique solution $\mu(\cdot)\in C^0(\R;\mathcal{P}_c(\R^d))$, where $\mathcal{P}_c(\R^d)$ is endowed with the weak-$*$ topology of measures. 
Moreover, $t\mapsto \mu(t)$ is Lipschitz in the sense of the Wasserstein distance $W_p$. 
Moreover, if $\mu^0\in\mathcal{P}^{ac}_c(\R^d)$, then $\mu(t)\in\mathcal{P}^{ac}_c(\R^d)$ for every $t\in\R$.

Furthermore, for every $T>0$, there exists $C_T>0$ such that 
\begin{equation}\label{e-stabilitypde}
W_p(\mu(t),\nu(t))\leq e^{C_T t}W_p(\mu(0),\nu(0)),
\end{equation}
for all solutions $\mu$ and $\nu$ of \eqref{e-cauchy} in $C^0([0,T];\mathcal{P}_c(\R^d))$. 

Finally, the solution $\mu$ of the Cauchy problem \eqref{e-cauchy} can be made explicit as follows.
Let $\Phi(t)$ be the flow of diffeomorphims of $\R^d$ generated by the time-dependent vector field $f[\mu]$, defined as the unique solution of the Cauchy problem $\dot\Phi(t) = f[\mu(t)] \circ \Phi(t)$, $\Phi(0)=\mathrm{Id}_{\R^d}$, or in other words,
$$
\partial_t\Phi(t,x) = f[\mu(t)](t,\Phi(t,x)),\quad \Phi(0,x)=x.
$$
Then, we have
$$
\mu(t) = \Phi(t)\#\mu_0,
$$
that is, $\mu(t)$ is the push-forward of $\mu_0$ under $\Phi(t)$.
\end{theorem}

Theorem \ref{t-esistenza} can be generalized to mass-varying transport PDEs, that is, in presence of sources (see \cite{genwass}).
We now observe that Theorem \ref{t-esistenza} can be applied to \r{e-PDE} as well, under Assumptions {\bf (H)} and provided that the control $u$ be a Lipschitz function of the space variable for all times.

\begin{corollary} 
Under Assumptions {\bf (H)}, if $u$ is a uniformly Lipschitz function of $x$ on the time interval $[0,\theta]$, satisfying the constraint ${\bf (U)}$, then, given any initial data $\mu(0)=\mu_0\in\PR$, the equation \r{e-PDE} has a unique solution $\mu(\cdot) \in C^0([0,\theta],\PR)$.
Moreover, if $\mu_0\in\PRac$, then $\mu(t)\in\PRac$ for every $t\in[0,\theta]$.
Denoting by $\Psi$ the flow of diffeomorphims of $\R^d$ generated by the time-dependent vector field $f[\mu]+u(t,x) g[\mu]$, we have $\mu(t) = \Psi(t)\#\mu_0$.
\end{corollary}

\bproof 
It suffices to check that Theorem~\ref{t-esistenza} can be applied to the vector field $f[\mu]+u(t,x) g[\mu]$. As already stated, the existence of a uniform bound $M$ for $\|f+ug\|_{L^\infty}$ is a consequence of the uniform Lipschitz property and of the uniform boundedness of the support of both $f[\mu]$ and $g[\mu]$, together with the bound $\|u\|_{L^\infty}\leq 1$ imposed by {\bf (U)} in~\eqref{e-strength}.
Similarly, we have a uniform bound on the Lipschitz constant $\Lip(f+ug)$. Indeed, by~\eqref{e-spacecon},
\begin{align}
 \mathrm{Lip}(f[\mu(t)]+ug[\mu(t)])&\leq L+\mathrm{Lip}_x(u)\|g[\mu(t)]\|_{L^\infty}+\|u\|_{L^\infty} \mathrm{Lip}(g[\mu(t)]) \nonumber \\
&\leq 2L+\mathrm{Lip}_x(u) M.
\label{e-stimaLip}
\end{align}
Finally, we have
\begin{align*}
\|f[\mu]+u g[\mu]-f[\nu]+u g[\nu]\|_{L^\infty}&\leq \|f[\mu]-f[\nu]\|_{L^\infty}+ \|u\|_{L^\infty} \|g[\mu]-g[\nu]\|_{L^\infty}\\
&\leq W_p(\mu,\nu)+1\cdot QW_p(\mu,\nu)\\
&=2QW_p(\mu,\nu).
\end{align*}
This proves the corollary.
\eproof

We end this section with an estimate of the $L^\infty$ norm of the solution $\mu(t)$ to \r{e-pdebase}, when it is absolutely continuous with respect to the Lebesgue measure.

\bp \label{p-Linfgenerale}
Let $\mu(\cdot)$ be the unique solution of \r{e-pdebase} for a given Lipschitz vector field $f$ with $\mu_0\in\Pac$. Then
\bqn
\frac{d}{dt}\|\mu(t)\|_{L^\infty}\leq \|\mu(t)\|_{L^\infty}\|\nabla\cdot f\|_{L^\infty}.
\eqnl{e-Linfgenerale}
\ep

\bproof 
The proof follows \cite[Proposition~3.1]{ha-tad}. Let $\rho(t)$ be the density of $\mu(t)$ with respect to the Lebesgue measure. For each $p\in[1,+\infty)$, by dropping the dependence with respect to time, we write 
$$
\frac{d}{dt}\int \rho^p\,dx=-p\int \rho^{p-1}\,\nabla\cdot (f\rho)\,dx=-p\int \Pt{ \rho^p \,\nabla\cdot f +\rho^{p-1} \Pa{f,\nabla \rho}}\,dx.
$$
Since $\nabla\cdot(f \rho^p)=\rho^p\,  \nabla\cdot f+\Pa{f,\nabla (\rho^p)}=\rho^p \nabla\cdot f+ p \rho^{p-1} \Pa{f,\nabla \rho}$, we infer that
\bqn
\frac{d}{dt}\int \rho^p\,dx&=&-\int(p-1)\rho^p \,\nabla\cdot f\,dx-\int \nabla\cdot(f\rho^p)\,dx.
\eqnn
The last term is zero as a consequence of the divergence theorem. Then 
$\frac{d}{dt}\|\rho\|_{L^p}^p\leq (p-1)\|\rho\|_{L^p}^p \|\nabla\cdot f\|_{L^\infty}$,
which in turn implies \r{e-Linfgenerale} as $p \to +\infty$.
\eproof

\section{Steepest descent under population constraint induces mass concentration}
\label{s-concentration}
In this section, we discuss a remarkable phenomenon for controlled equations of the form \r{e-PDE}: starting from a measure $\mu_0\in\Pac$, i.e., a measure that is absolutely continuous with respect to the Lebesgue measure in $\R^d$, a time-dependent choice of the control might drive the measure outside $\Pac$ in finite time, in particular with emergence of Dirac deltas. In fact, we will show that {\em such a phenomenon arises when trying to minimize a Lyapunov function} $V$, in particular when one chooses the control $u(t)$ as the instantaneous minimizer of  the Lie derivative of $V$ as time evolves. This example also shows that some key ideas coming from control of finite-dimensional systems cannot be extended straightforwardly to infinite dimension.

In this section, we discuss the interest and the drawbacks of a control constraint different than {\bf (U)}, namely the following:

\noindent\fbox{%
    \parbox{0.99\textwidth}{%
    \begin{center}
    {\bf Alternative control Constraints (U')}
    \end{center}

Fix $c>0$. For each time $t\geq0$ it holds:
\bqn
&&\mbox{{\bf Sparsity population constraint:}}\hspace{2cm} \int_{\omega(t)} d\mu(t) \leq c,\label{e-popcon}\\
&&\mbox{{\bf Finite strength:}} \hspace{47mm}  \| u(t,.)\|_{L^\infty}\leq 1.
\eqn

}}

\vspace{3mm}

The population constraint represents the idea of acting on a small part of the crowd itself, and not on a small part of the configuration space, as we require in the space constraint in {\bf (U)}. Even though the sparse population constraint is interesting from the theoretical point of view, it has a surprising drawback on the modeling point of view: when a crowd is extremely concentrated, the constraint {\bf (U')} implies that the control cannot act on the crowd anymore. This is somehow unnatural, since a crowd that is already concentrated is the best configuration to steer. On the other hand, the space constraint {\bf (U)} permits to act on the whole crowd, when it is concentrated in a set of size $c$, i.e., exactly when it is concentrated.

We now show that the population constraint also induces some formal mathematical problems when using a sparse Jurdjevic--Quinn approach. Consider the following system on the real line:
\bqn
\partial_t\mu+\partial_x( u \mu)=0.
\eqnl{e-example}
This is a particular case of \r{e-PDE} with $f=0$ and $g=1$. We consider the initial data $\mu_0=\chi_{[0,1]}$, i.e., a uniform probability density on the interval $[0,1]$. We consider the Lyapunov function 
\bqn
V[\mu]=\int x^2 d\mu(x),
\eqnn
i.e., the second moment with respect to zero. We have $\L_{f} V=0$, and we have $\L_{u g}V[\mu]=0$ for $\mu=\delta_0$ only, i.e., $\mathcal{Z}=\Pg{\delta_0}$.  Then, minimizing $V$ is equivalent to steer $\mu(t)$ to the Dirac mass $\delta_0$.

We now apply a rough form of the steepest descent method to the problem of minimizing $V$: given the initial measure $\mu_0$, we look for a control function $u$ that maximizes the descent $\L_{u g}V[\mu]$, while taking into account the control constraints {\bf (U')}. An easy computation shows that no optimal choice for $u$ exists. Indeed, for every $\eps>0$, consider the $C^\infty$ function 
$$
u_\eps(x)=
\begin{cases}
-1 & \mbox{ for } x\in [1-c+\eps,1],\\
0 & \mbox{ for } x\in (-\infty, 1-c]\cup [1+\eps,+\infty),\\
C^\infty-\mbox{spline with values in } [-1,0] & \mbox{ for } x\in [1-c,1-c+\eps] \cup [1,1+\eps].
\end{cases}
$$
Then, for a sufficiently small time $t>0$, each particle $x\in(1-c+\eps,1]$ is displaced to $x-t$ while each particle $x\in[0,1-c]$ undergoes no displacement. The particles in the small interval $[1-c,1-c+\eps]$ are displaced toward $1-c$, then giving a reduction of the value of the functional $V$. Then, we have 
\begin{align*}
 \L_{u_\eps g}V[\mu]= &\ 
 \frac{d}{dt}_{|_{t=0}} \int_{1-c}^{1-c+\eps}  (x+tu(x)+\mathrm{o}(t))^2 \, d\mu(x) +\frac{d}{dt}_{|_{t=0}} \int_{1-c+\eps}^1 (x-t)^2 \, d\mu(x) \\
 \leq &  -2 \int_{1-c+\eps}^1 x \, d\mu(x).
\end{align*}
As a consequence, by decreasing the parameter $\eps>0$, one can reach a larger decrease of $V$. Nevertheless, the limit for $\eps \to 0$ would result in the discontinuous control function $u_0=\chi_{[1-c,1]}$, for which the solution to the corresponding dynamics \r{e-example} does not satisfy Assumptions\footnote{Existence and uniqueness for the solution of \r{e-PDE} in small times, with possibly discontinuous controls, can be derived from results of \cite{ambrosio,paola}. Nevertheless, it is shown in \cite{paola} that, under such assumptions, one can have formation of singularities such as Dirac deltas in finite time.} {\bf (H)}.

The fact that a maximizer of the steepest descent does not exist in the space of Lipschitz functions can be  overcome by fixing a value $\eps_0>0$ and applying the control $u_{\eps_0}$ over a small interval of time $[0,t_0]$ with $t_0<\eps_0$. As a result, the component of the measure $\mu$ with $x\in [1-c,1]$ concentrates in the interval $[1-c,1-t_0]$, while its density keeps being constantly equal to 1 for $x\in [0,1-c]$. At time $t$, one can observe that the largest contribution to $V$ comes anyway from the mass in the interval $[1-c,1-t_0]$, on which the control already acted. Moreover, the mass in such interval keeps being $c$.

For the reasons described above, any strategy maximizing the descent of $V$ acts on the mass in the interval $[1-c,1]$ for all times. In particular, we can choose a sequence $\eps_i$ acting on the time interval $[t_{i-1},t_i]$, with the condition $\eps_i<c-t_i$ for all $t\in[0,c)$. Applying the time-dependent control $u(t,x)=u_{\eps_i}(x)$ for $t\in[t_{i-1},t_i)$, the corresponding solution $\mu(t)$ of \r{e-example} has support in $[0,1-t]$. More precisely, it has the following structure: the measure keeps having density 1 in the interval $[0,1-c]$, while the rest of the mass $c$ is contained in the interval $[1-c,1-t]$.

Then, the solution $\mu(t)$ converges as $t\nearrow c$ to the singular measure 
\bqn
\mu(c)=\chi_{[0,1-c]} +c\delta_{1-c}.
\eqnl{e-singular} 
This is  not in contradiction with the fact that any solution of \r{e-pdebase} with initial data $\mu_0\in\Pac$ and Lipschitz vector field $f$ satisfies $\mu(t)\in\Pac$. Indeed, the condition $\eps(t)<c-t$ implies $\lim_{t\to c} \eps(t)=0$, hence the control $u(t,x)$ converges to a non-Lipschitz function.

Starting from the singular measure \r{e-singular} at time $t=c$, one finds several problems to steer it towards the minimizer $\delta_0$ of $V$. First, the main contribution to $V$ is given by the Dirac delta $c\,\delta_{1-c}$: then, the control set $\omega$ chosen to maximize the descent for $V$ would certainly contain such mass. But the condition $1-c\in \omega$ together with $u(1-c)\neq 0$ would directly impose to choose $\omega$ containing a whole neighbourhood of $1-c$. This would imply $\int_\omega \mu(c)>c$, hence the population constraint would be automatically violated.

This would in turn enforce us to focus our control on the absolutely continuous part, possibly leading to the formation of a new Dirac delta $c\,\delta_{1-2c}$, and so on. The final result would be a set of Dirac deltas, not concentrated at $0$, on which control with population constraint cannot be applied.

Summing up, the steepest descent method with population constraint in {\bf (U')} might not steer the measure to a configuration in $\mathcal{Z}$, but rather to a configuration in which the population constraint itself may not be satisfied.

\section{Proof of Theorem~\ref{t-main}}
\label{s-proof}
For the moment, we assume that the solution $\mu(\cdot)$ of \eqref{e-PDEt}, with the control strategy defined by Theorem \ref{t-main}, is well defined on $[0,\theta]$ with $\theta\in(0,+\infty]$, and we establish some lemmas describing its evolution. Recall that $\supp(\mu_0)\subset \overline{B(0,R)}$.

\bl\label{l-compact} 
We have $\supp(\mu(t))\subset \overline{B(0,R)}$ for every $t\in[0,\theta]$.
\el

\bproof 
Since the vector field $f+ug$ is zero outside $\overline{B(0,R)}$, the corresponding flow $\phi^t(\cdot)$ coincides with the identity in $\R^d\setminus\overline{B(0,R)}$. Since we have $\mu(t)=\phi^t_\#\mu_0$ by Theorem \ref{t-esistenza}, we get that, for any Borel set $E$ satisfying $E\cap \overline{B(0,R)}=\emptyset$, we have $\mu(t)(E)=\mu_0(\phi^{-t}(E))=\mu_0(E)=0$. The lemma follows.
\eproof

Recall that $L$ is the Lipschitz constant for $f[\mu]$ and $g[\mu]$ given in~\eqref{e-regPDE}, and recall that $\|f[\mu]\|_{L^\infty} \leq M$ and $\|g[\mu]\|_{L^\infty}\leq M$ for every $\mu \in \PRac$ (as in Remark~\ref{rk:bound}).

\bl\label{l-Linf} 
We have $\|\mu(t)\|_{L^\infty}\leq e^{d\theta(2L+M\theta)}\|\mu(0)\|_{L^\infty}$, for every $t\in[0,\theta]$.
\el

\bproof 
Since the vector field $f+ug$ satisfies \r{e-stimaLip}, and since $\Lip_{x}(u(t,\cdot))\leq \frac{1}\eta\leq t\leq \theta$, the lemma follows from Proposition \ref{p-Linfgenerale}.
\eproof

\bl\label{lem:cont} 
The function $(t,a,b,\eta)\mapsto s_t(a,b,\eta)$ is continuous with respect to $t$, and uniformly Lipschitz with respect to $(a,b,\eta)$ on $\cup \{ \Omega_t, t\in[0,\theta]\}$.

\el
\bproof 
Let us first establish the Lipschitz property for $(a,b,\eta)\in \cup\{\Omega_t, t\in[t_n,t_{n+1}]\cap[0,\theta]\}$. Note that the condition $|\omega(a,b,\eta)|\leq c$ implies that $\eta\leq \frac{c}2$. Besides, we have $\eta\geq t^{-1}\geq \theta^{-1}$.
By definition of $U(a,b,\eta)$ in \r{e-U}, with simple geometric arguments, it is clear that in the 1D case we have
\bqn
\begin{split}
\| U(a,b,\eta)-U(a',b',\eta')\|_{L^1(dx)} &\leq 
\|U(a,b,\eta)-U(a,b,\eta')\|_{L^1(dx)}+\|U(a,b,\eta')-U(a',b',\eta')\|_{L^1(dx)} \\&\leq
|\eta-\eta'|+|a-a'|+|b-b'|,
\end{split}
\eqnn
where $dx$ is the standard Lebesgue measure on $\R$. By applying the estimate componentwise, the same result follows in dimension $d$. By Lemma~\ref{l-Linf}, there exists $P>0$ such that $\|\mu(t)\|_{L^\infty}\leq M$ for every $t\in[0,\theta]$, hence
$$
\| U(a,b,\eta)-U(a',b',\eta')\|_{L^1(\mu(t))}\leq P(|a-a'|+|b-b'|+|\eta-\eta'|),
$$
and thus \r{e-diff} implies that $s_t(a,b,\eta)$ is Lipschitz with respect to $(a,b,\eta)$, with Lipschitz constant $KP$. 

The function  $s_t(a,b,\eta)$ is continuous with respect to $t$, as a consequence of the continuity of $\L_{ug}V$ given by \r{e-contf}.
\eproof

\bl \label{l-max} 
Assume that $\mu(t_n)\in \Pac$. If $\Omega_{t_n}$ is nonempty then $s_{t_n}(a,b,\eta)$ has a maximizer in $\Omega_{t_n}$.
\el
\bproof It suffices to observe that $\Omega_{t_n}$ can be considered as being compact: indeed, each choice $(a,b,\eta)\in\Omega_{t_n}$ can be replaced by an equivalent choice $(a',b',\eta)$ with $a,b\in B_0(R+2c)$ since  $\|U(a,b,\eta)-U(a',b',\eta)\|_{L^1(\mu)}=0$, since $\mu$ has zero mass outside of $\overline{B(0,R)}$. In other terms, one can restrict the choice of the parameters $a,b$ to a compact set. Similarly, we have $\eta\in[t_n^{-1},\frac{c}2]$. Since, by Lemma~\ref{lem:cont}, $s_{t_n}$ is a continuous function of its arguments, then it admits a maximizer.
\eproof

We are now in a position to prove Theorem \ref{t-main}. 
We split the proof into three steps:
\begin{description}
\item[Step 1.] For each time $t_n$, the $n$-th step of the algorithm univocally determines a control satisfying the constraint {\bf (U)}, the corresponding solution of \r{e-PDE} and a time $t_{n+1}>t_n$.
\item[Step 2.] We have $t_n\rightarrow+\infty$.
\item[Step 3.] We have $\lim_{t\to+\infty}\L_{f+ug}V[\mu(t)]=0$. This fact, together with the choice of maximizing controls and of the hysteresis, provides convergence to the sets in which the maximizers of $s_t(a,b,\eta)$ give zero control. Since the constraint $\eta\geq t^{-1}$ is negligible for $t\to+\infty$, the strategy provides convergence of $\mu(t)$ to  $\mathcal{Z}$.
\end{description}

\medskip

\noindent {\sc Proof of Step 1.} Let us prove that the algorithm of Theorem \ref{t-main} univocally defines a control strategy, by induction. We have $t_0=0$ and $\mu(0)=\mu_0\in\Pac$. Let us prove that, for a given time $t_n$, the time $t_{n+1}$ is well defined and satisfies $t_{n+1}>t_n$.

We first observe that the control $\chi_\omega u$ is a well-defined function, Lipschitz with respect to $x$. Setting $\tilde t=\frac12 \Pt{\frac{c}{|B_0(1)|}}^{1/d}$, we note that $\Omega_t=\emptyset$ for every $t\in[0,\tilde{t})$ since any function of the form $\chi^\eta_{[a,b]}$ has a support of size larger than $(2\eta)^d |B_0(1)|$. For $t\geq \tilde{t}$, the set $\Omega_t$ is nonempty and Lemma~\ref{l-max} yields the existence of a maximizer $(a^*,b^*,\eta^*)$ in $\Omega_t$.
We thus have two cases:
\bi
\i If $s_{t_n}(a^*,b^*,\eta^*)< t_n^{-1}$ or $\Omega_t$ empty, then the control $\chi_\omega u=0$ is well defined and is Lipschitz.
\i If $s_{t_n}(a^*,b^*,\eta^*)\geq t_n^{-1}$, then the control $\chi_\omega u=U(a^*,b^*,\eta^*)$ is well defined and is Lipschitz, as a consequence of the Lipschitz property in~\eqref{e-chi}.
\ei
Let us now prove that there exists a unique minimum $t_{n+1}$ defined by the algorithm, and that it satisfies $t_{n+1}>t_n$. For $t\in [0,\tilde{t})$, there is nothing to prove, since $t_1\geq \tilde{t}$. For $t\geq \tilde{t}$, we have two cases:
\bi
\i If $s_{t_n}(a^*,b^*,\eta^*)< t_n^{-1}$, then the set 
$$
A=\Pg{t\geq t_n \ \mid \ s_t(a,b,\eta)\geq t^{-1}\mbox{~~for some~~} (a,b,\eta)\in\Omega'_t}
$$ 
is closed or empty. If it is nonempty, there exists a minimal element $t_{n+1}\geq t_n$. Moreover, $t_{n+1}\neq t_n$, since $s_{t_n}(a,b,\eta)\leq s_{t_n}(a^*,b^*,\eta^*)< t_n^{-1}$ for all $(a,b,\eta)\in\Omega_{t_n}$. If $A$ is empty then $+\infty=t_{n+1}>t_n$.

\i If $s_{t_n}(a^*,b^*,\eta^*)\geq t_n^{-1}$, then, similarly to the previous case, since the function $s_t(a^*,b^*,\eta^*)$ is continuous with respect to the time $t$, the set 
$$
A'=\Pg{t\geq t_n \ \mid \ s_{t}(a^*,b^*,\eta^*)\leq \frac{t_n^{-1}}2 }
$$ 
is closed or empty, and it does not contain $t_n$.
We now consider the set 
$$B'=\Pg{\ba{l}t\geq t_n \mbox{~such that~there exists~}(\bar a,\bar b,\bar \eta)\in \Omega'_t\mbox{~~for which~}\\
s_t(\bar a,\bar b,\bar \eta)\geq (1-h)^{-1} s_t(a^*,b^*,\eta^*)\geq (1-h)^{-1}\frac{t^{-1}}2.\ea}. 
$$
Let us prove that it is closed and that $t_n\not\in B'$. Take a sequence $(t^i,a^i,b^i,\eta^i)$ such that $t^i\in B'$ is a sequence converging to some $\hat t$, and $(a^i,b^i,\eta^i)$ satisfy $s_{t^i}(a^i, b^i,\eta^i)\geq (1-h)^{-1}s_{t^i}(a^*,b^*,\eta^*)\geq \frac{\tau(t^i)}2$. 
Observing that the compact set $\Omega'_t$ varies smoothly with respect to time, we can restrict ourselves to a sequence $(a^i,b^i,\eta^i)$ converging to some $(\bar a,\bar b,\bar \eta)\in\Omega'_{\bar t}$. Then, by continuity of $s_t$, we have $s_{\bar t}(\bar a,\bar b,\bar \eta)\geq (1-h)^{-1}s_{\bar t}(a^*,b^*,\eta^*)$, hence $\bar t\in B'$. Moreover, $\bar t\neq t_n$, otherwise 
$(a^*,b^*,\eta^*)$ would not be a maximizer of $s_{t_n}$.

Since both $A'$ and $B'$ are closed or empty, not containing $t_n$, then $A'\cup B'$ is closed or empty and does not contain $t_n$. If it is closed, then it admits a minimal element $t_{n+1}>t_n$; if it is empty, then we have $+\infty=t_{n+1}>t_n$. 
\ei

\medskip

\noindent {\sc Proof of Step 2.} We now prove that the sequence $t_n$ of times given by the algorithm converges to $+\infty$. Since $t_n$ is increasing, it has a limit $T$. By contradiction, if $T<+\infty$, then $\mu(t)$ is defined for every $t\in[0,T]$. Indeed, since $\|f+ug\|_{L^\infty}\leq 2M$, the curve $t\mapsto\mu(t)$ is Lipschitz, and thus $\mu(T)$ is well defined.

If we have $s_{t_n}(a^*,b^*,\eta^*)<t_n^{-1}$ at time $t_n$, then at the next time $t_{n+1}$ we must have $s_{t_{n+1}}(\bar a, \bar b,\bar \eta)\geq t_{n+1}^{-1}$ for some $(\bar a, \bar b,\bar \eta)\in\Omega_{t_{n+1}}$, by definition of the algorithm itself. As a consequence, the sequence $t_n$ converging to $T$ contains an infinite number of times $t_{n_i}$ such that $s_{t_{n_i}}(a^i,b^i,\eta^i)\geq t_{n_i}^{-1}$, where $(a^i,b^i,\eta^i)$ is a maximizer of $s_{t_{n_i}}$ in $\Omega_{t_{n_i}}$.

The sequence $(a^i,b^i,\eta^i)$ is bounded, and its converging subsequences have their limit in $\Omega_T$. Indeed, we can restrict ourselves to $(a^i,b^i)\in B_0(R+2c)$, and we have $\eta^i\geq T^{-1}$ and $\eta^i\leq \frac{c}2$. Hence, taking a subsequence if necessary, we have the existence of a limit $(\hat a, \hat b,\hat \eta)\in\Omega_T$.

Observe now that, at time $t_{n_i+1}$, one of the two conditions leading to switching of the control holds. Since the sequence $t_{n_i+1}$ has an infinite number of terms, at least one of the conditions holds for an infinite subsequence (that we do not relabel). We show now that this is in contradiction with the fact that $t_n$ converges to a finite time $T$. We have two cases:
\bi
\i If $s_{t_{n_i}}(a^i,b^i,\eta^i)\geq t_{n_i}^{-1}$ and $s_{t_{n_i+1}}(a^i,b^i,\eta^i)\leq \frac{t_{n_i}^{-1}}2$, then, taking a subsequence converging to $(\hat a, \hat b, \hat \eta)\in\Omega_T$, we have a contradiction with the continuity of $s_T$ in $(\hat a, \hat b, \hat \eta)$. Indeed, we have
\bqn
s_T(\hat a, \hat b, \hat \eta)=\lim_{i\to+\infty} s_{t_{n_i}}(a^i,b^i,\eta^i)\geq t_{n_i}^{-1}>\frac{t_{n_i}^{-1}}2\geq \lim_{i\to+\infty} s_{t_{n_i}+1}(a^i,b^i,\eta^i)=s_T(\hat a, \hat b, \hat \eta).
\eqnn
\i If there exists $(\bar a^i,\bar b^i,\bar \eta^i)\in \Omega'_{t_{n_i+1}}$ such that 
\bqn
s_{t_{n_i+1}}(a^i,b^i,\eta^i)\leq (1-h) s_{t_{n_i+1}} (\bar a^i,\bar b^i,\bar \eta^i),
\eqnl{e-min1}
then, for $n_i\to+\infty$, we have $\Omega'_{t_{n_i+1}}\subset \Omega_{t_{n_i}}$, since $2t_{n_i+1}^{-1}\geq t_{n_i}^{-1}$, as a consequence of the fact that $\lim_{i\to+\infty}(t_{n_i+1}-t_{n_i})=0$. Since $(a^i,b^i,\eta^i)$ is a maximizer of $s_{t_{n_i}}$ in $\Omega_{t_{n_i+1}}$, we have
\bqn
s_{t_{n_i}}(\bar a^i,\bar b^i,\bar \eta^i)\leq 
s_{t_{n_i}}(a^i,b^i,\eta^i).
\eqnl{e-min2}
One can take a converging subsequence of $(\bar a^i,\bar b^i,\bar \eta^i)$, for the same reasons given above for the sequence $(a^i,b^i,\eta^i)$. Denoting by $(\bar a^*,\bar b^*,\bar \eta^*)$ and $(a^*,b^*,\eta^*)$ the two limits, and using continuity of $s_t(a,b,\eta)$ with respect to all its arguments, we get from \r{e-min1}-\r{e-min2} that
\bqn
s_{T}(a^*,b^*,\eta^*)
\leq (1-h) s_{T} (\bar a^*,\bar b^*,\bar \eta^*)
\leq (1-h) s_{T}(a^*,b^*,\eta^*),
\eqnn
which is in contradiction with $s_{T}(a^*,b^*,\eta^*)\geq T^{-1}>0$.
\ei
Then $t_n$ cannot converge to a finite value $T$. Therefore either $t_n\to+\infty$ or there exists a $t_n$ such that $t_{n+1}=+\infty$. In both cases, the control strategy is defined for every $t\in[0,+\infty)$.

\medskip

\noindent {\sc Proof of Step 3.} 
It remains to prove that $\mu(t)$ converges to $\mathcal{Z}$. This is the hardest part of the proof, in which the choice of the admissible controls in $\Omega_t$ plays a crucial role.

Thanks to Step 2, we have, for every time, $\mu(t)\in \PR$, that is compact with respect to the weak topology, which coincides with the topology of the Wasserstein distance. Then Assumptions {\bf (H)} imply that $V$ is a continuous function and thus is bounded below.

We now prove that the function $V(t)= V[\mu(t)]$ is differentiable for almost every $t$, and that it satisfies 
$\dot V(t) = \lim_{t\to 0}\frac{V[e^{t(f+ug)}\mu]-V[\mu]}{t} \leq 0$. Differentiability on the open time interval $(t_n,t_{n+1})$ follows from the fact that $\dot V(t)=\L_{f+ug}V[\mu(t)]$ is continuous, as a consequence of Assumptions {\bf (H)}. Clearly, the set of times $t_n$ on which differentiability is not ensured is countable, hence $V(t)$ is differentiable for almost every $t$.

For $t\in(t_n,t_{n+1})$, we have $\dot V=\L_{f+ug}V[\mu]=\L_fV[\mu]+\L_{ug}V[\mu]$. If at time $t_n$ the algorithm defines the control $\chi_{\omega}u\equiv 0$, then clearly 
\bqn
\dot V=\L_f V[\mu]\leq 0 ,
\eqnl{e-dotV1} 
for every $t\in(t_n,t_{n+1})$. If instead the control given by the algorithm is $\chi_\omega u$ in \r{e-control}, we have 
\bqn
\dot V=\L_f V[\mu]+\L_{\chi_\omega u g} V[\mu]\leq -\mathrm{sign}(\L_{U(a^*,b^*,\eta^*)g[\mu(t_n)]}V[\mu_{t_n}]) (\L_{U(a^*,b^*,\eta^*)g[\mu(t)]}V[\mu_{t}]).
\eqnl{e-dotV2}
It is clear that, at the beginning of the interval, we have
\bqn
\lim_{t\to t_n^+}\mathrm{sign}(\L_{U(a^*,b^*,\eta^*)g[\mu(t_n)]}V[\mu_{t_n}]) (\L_{U(a^*,b^*,\eta^*)g[\mu(t_n)]}V[\mu_{t_n}])\geq s_{t_n}(a^*,b^*,\eta^*)\geq \tau(t_n)>0,
\eqnn
and hence $\dot V(t_n^+)<0$. Since $\dot V(t)$ is a continuous function, we either have $\dot V(t)<0$ for every $t \in (t_n,t_{n+1})$, or there exists $t\in (t_n,t_{n+1})$ such that $\L_{U(a^*,b^*,\eta^*)g[\mu(t)]}V[\mu(t)]=0$. This is equivalent to state that $s_t(a^*,b^*,\eta^*)=0$, which is in contradiction with $s_t(a^*,b^*,\eta^*)>\frac{\tau(t)}2>0$ for every $t\in (t_n,t_{n+1})$, by definition of the time $t_{n+1}$.

We now prove that $\lim_{t\to \infty}\mu(t)\in \mathcal{Z}$. Since $\PR$ is compact, all sequences have limits. Consider a sequence $t^j\to\infty$ such that $\lim_{j\to\infty}\mu(t^j)=\mu^*$. We are going to prove that $\mu^*\in \mathcal{Z}$.

Since $V$ is continuous, bounded below and $\dot V(t)\leq 0$ for almost every $t$, we have $\lim_{t\to +\infty} V(t)=V^*$ for some $V^*$. The existence and continuity of the second-order derivatives $\L_{f+ug}\L_{f+ug}V[\mu]$ on the compact space $\PR$ implies the existence of a uniform bound on $\ddot V$. As a consequence, we have $\lim_{t\to\infty}\dot V(t)=0$. Since $\dot V\leq \L_f V\leq 0$ by either \r{e-dotV1} or \r{e-dotV2}, this in turn implies $\lim_{j\to\infty} \L_f V[\mu(t^j)]=0$, hence $\L_fV[\mu^*]=0$ by continuity of $\L_fV$.

We now prove that $\L_{ug}V[\mu^*]=0$ for all $u\in \U$. By contradiction, assume that there exists $u^*\in \U$ such that $|\L_{u^*g}V[\mu^*]|\neq 0$. Without loss of generality, by using \r{e-mult}, we assume that $\|u^*\|_{L^\infty}=1$. Similarly, by decomposing $u^*=u^+-u^-$ with $u^+,u^-$ non-negative Lipschitz functions, and using additivity of the Lie derivative, we can replace $u^*$ with either $u^+$ or $u^-$ and assume that it is nonnegative and that $|\L_{u^*g}V[\mu^*]|= C^*\neq 0$. Finally, by observing that $\mu^*$ has compact support, we can replace $u^*$ with a nonnegative Lipschitz function with compact support.

We now approximate $u^*$ in $L^\infty$ by a family of functions of the form $\sum_{i=1}^I k^i\chi^{\eta}_{[a^i,b^i]}$, where the number $I$ of terms depends on the approximation error, but not on the (sufficiently small) parameter $\eta$. 

For simplicity, we only give the construction in the 1D case. 

Since $u^*$ is Lipschitz with bounded support, it is Riemann integrable. In particular, by using an approximation of $u^*$ from below, we have the following: take a grid step $\dx$ and define  the  rectangles $k^i\chi_{[\tilde a^i,\tilde b^i]}$ with $\tilde b^i-\tilde a^i=\dx$, for which
\bqn
\sum_{i=1}^I k^i\chi_{[\tilde a^i,\tilde b^i]}\leq u^*\mbox{~~and~~}\|u^*-\sum_{i=1}^I k^i\chi_{[\tilde a^i,\tilde b^i]}\|_{L^\infty}\leq \eps.\label{e-Riemann}
\eqn
for some $\eps$. The Riemann integrability of $u^*$ implies that, for any $\eps>0$ there exists $\dx$ such that \r{e-Riemann} is satisfied. Note that $\|u^*\|_{L^\infty}=1$ also implies $k^i\leq 1$.

We now prove that we can replace $k^i\chi_{[\tilde a^i,\tilde b^i]}$ with their mollified version $k^i\chi^\eta_{[a^i,b^i]}$ for any sufficiently small $\eta$, while keeping \r{e-Riemann} satisfied. We provide here the explicit construction. First denote with $L'$ the Lipschitz constant of $u^*$. To replace $k^i\chi_{[\tilde a^i,\tilde b^i]}$ with $k^i\chi^\eta_{[a^i,b^i]}$, we have two cases:
\bi
\i If $\tilde b^{i-1}<\tilde a^i$, then keep both $b^{i-1}=\tilde b^{i-1}$ and $a^i=\tilde a^i$, and choose $\eta\leq \min\Pg{\frac{\tilde a^i-\tilde b^{i-1}}2,\frac{k^{i-1}}{L'},\frac{k^{i}}{L'}}$.
\i If $\tilde b^{i-1}=\tilde a^i$, then choose $\eta\leq\min\Pg{\frac{k^{i-1}}{L'},\frac{k^{i}}{L'}}$. If $k^{i-1}>k^{i}$, then define $b^{i-1}=\tilde b^i$ and $a^i=\tilde a^i+\eta$. Otherwise, take $b^{i-1}=\tilde b^i-\eta$ and $a^i=\tilde a^i$.
\ei

\begin{figure}[ht]
\begin{center}
\includegraphics[width=13cm]{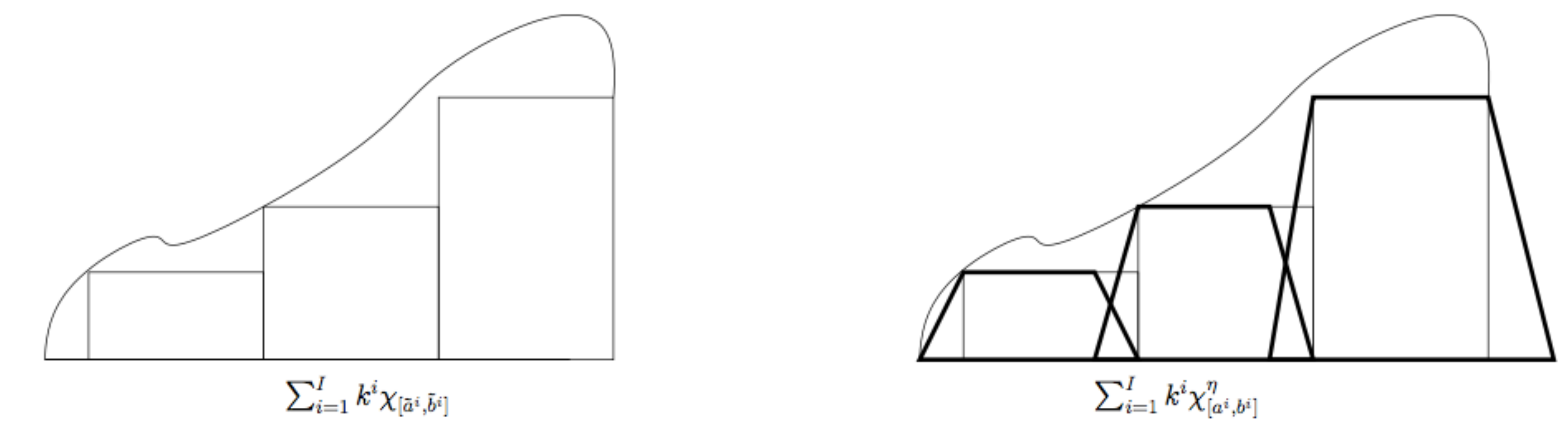}
\caption{Construction of the approximation $k^i\chi^\eta_{[a^i,b^i]}$. }
\label{fig:chieta}
\end{center}
\end{figure}

Note that the constraints imposed on $\eta$ are higher bounds, and they are in finite number. Then, $\eta$ can be chosen in a whole interval $(0,\eta']$, where $\eta'$ depends on $\eps$ only. By construction, we have $\sum_{i=1}^I k^i\chi_{[\tilde a^i,\tilde b^i]}\leq \sum_{i=1}^I k^i\chi^\eta_{[ a^i, b^i]}\leq u^*$, thus both conditions in \r{e-Riemann} are satisfied. Note that this construction depends only on  $u^*$ and not on $\mu^*$. We denote this new function by $u^\eta=\sum_{i=1}^I k^i\chi^\eta_{[ a^i, b^i]}$. 

Consider now the sequence of measures $\mu(t^j)$ converging to $\mu^*$. If $t^j$ is a switching time $t_n$ for the algorithm, replace it with a slightly larger time. Then we can assume that  $\mu(t^j)$ keeps converging to $\mu^*$, with no switching times. By continuity of $\L_{u^*g}V[\mu]$ given by \r{e-contf}, we have $|\L_{u^*g}V[\mu(t^j)]|\geq \frac{C^*}2$ for sufficiently large indices. Note that 
\bqn
|\L_{(u^*-u^\eta)g}V[\mu(t^j)]|\leq K \|u-u^\eta\|_{L^1(\mu(t^j))}\leq K\|u-u^\eta\|_{L^\infty} \leq K \eps.
\eqnn
Then, choose a sufficiently small $\eps\leq \frac{C^*}{2K}$ and a corresponding $\dx>0$ such that \r{e-Riemann} is satisfied. If $\dx\geq \frac{c}2$, then replace it with $\frac{c}2$: by integrability of $u^*$, \r{e-Riemann} is still satisfied when refining the grid. Then, the previous construction shows that there exists $\eta'$ such that $u^\eta$ satisfies \r{e-Riemann} for every $\eta\in(0,\eta']$. Choose then $\eta^*\in(0,\eta']$. For a sufficiently large $j$, we have $
|\L_{u^{\eta^*} g}V[\mu(t^j)]|\geq \frac{C^*}2$.

Note now that $u^{\eta^*}$ is a sum of $I$ terms of the form $k^i\chi^{\eta^*}_{[ a^i, b^i]}$. Then, for each $j$ there exists at least one term such that 
\bqn
|\L_{\chi^{\eta^*}_{[ a^i, b^i]} g}V[\mu(t^j)]|\geq |\L_{k^i\chi^{\eta^*}_{[ a^i, b^i]} g}V[\mu(t^j)]|\geq \frac{C^*}{2I},
\eqnn
where we have used that $k^i\leq 1$ and that the number $I$ does not depend on the parameter $\eta^*$. Observe now that, for a sufficiently large $T$, we have $(a^i,b^i,\eta^*)\in\Omega'_t$ for every $t\geq T$. Similarly, taking a larger $T$ if necessary, we have $s_{t^j}(a^i,b^i,\eta^*)\geq  \frac{C^*}{2I}\geq t_j^{-1}$. As a consequence, the control algorithm provides a maximizer $(a^j,b^j,\eta^j)\in \Omega'_{t^j}$ of $s_{t^j}$, for which
\bqn
s_{t^j}(a^j,b^j,\eta^j)\geq (1-h)s_{t^j}(a^i,b^i,\eta^*)\geq (1-h)\frac{C^*}{2I}. 
\eqnn
The corresponding derivative satisfies 
\bqn
\dot V(t^j)=\L_{f+U(a^j,b^j,\eta^j)g}V[\mu(t^j)]\leq 0 - (1-h)\frac{C^*}{2I}. 
\eqnn
In particular, this is in contradiction with the fact that $\lim_{t\to\infty}\dot V(t)=0$.
The theorem is proved.

%
%

\section{Generalization to several controls}\label{sec:gen}
In this section, we show how to extend our result to a transport equation with a finite number of controlled vector fields,
$$
\partial_t \mu +\nabla\cdot \left( \left( f\Pq{\mu}+ \sum_{i=1}^m \chi_{\omega_i} u_i g_i[\mu] \right) \mu \right) =0.
$$
Now, Assumptions {\bf (H)} are done for all vector fields $f,g_1,\ldots,g_m$. In this setting, one can require the control constraint {\bf (C)}, together with the following additional constraint.\\

\noindent\fbox{%
    \parbox{0.99\textwidth}{%
    \begin{center}
    {\bf Componentwise sparsity constraint}
\end{center}
    
For every $t\in[0,+\infty)$, there exists at most one index $i$ such that $u_i(t,\cdot)$ is not identically zero.
}}

\vspace{2mm}

This sparsity constraint was first considered in the finite-dimensional setting for crowd models in \cite{CS-control,CS-control-1}. We recently generalized the Jurdjevic--Quinn stabilization method with this additional constraint in \cite{finiteJQ}.
Under such additional assumptions, we can adapt the control algorithm of Theorem \ref{t-main} as follows.\\

\noindent\fbox{%
    \parbox{0.99\textwidth}{%
    
     \begin{center}
    {\bf Step $n$}
\end{center}
    
    At time $t_n$, choose the maximizer $(a^*,b^*,\eta^*,i^*)$ of $|\L_{U(a,b,\eta)g_i}V[\mu(t_n)]|$ in the space $\Omega_{t_n}\times\Pg{1,\ldots,m}$. Observe that we only added the index $i$ in the maximization process.

Then, we have two cases:
\bi
\i If $|\L_{U(a^*,b^*,\eta^*)g_{i^*}}V[\mu(t_n)]|< t_n^{-1}$ or if $\Omega_{t_n}$ is empty, then choose the zero control 
$$
\chi_\omega u(t,x)\equiv 0
$$ 
and let the measure $\mu(t)$, starting at $\mu(t_n)$, evolve according to \r{e-PDEt} over the time interval $[t_n,t_{n+1}]$, where $t_{n+1}$ is the smallest time greater than $t_n$ such that there exists $(\bar a,\bar b,\bar \eta,\bar i)\in\Omega'_{t_{n+1}}\times\Pg{1,\ldots,m}$ for which $|\L_{U(\bar a,\bar b,\bar \eta)g_{\bar i}}V[\mu(t_{n+1})]| \geq 2t^{-1}$.
\i If $|\L_{U(a^*,b^*,\eta^*)g_{i^*}}V[\mu(t_n)]|\geq t_n^{-1}$, then choose the control 
$$
\chi_{\omega}u(t,\cdot)=-U(a^*,b^*,\eta^*)\ \mathrm{sign}(\L_{U(a^*,b^*,\eta^*)g_{i^*}[\mu(t)]} V[\mu(t)])
$$
and let the measure $\mu(t)$, starting at $\mu(t_n)$, evolve according to \r{e-PDEt} over the time interval $[t_n,t_{n+1}]$, where $t_{n+1}$ is the smallest time greater than $t_n$ satisfying at least one of the following conditions:
\bi
\i either $|\L_{U(a^*,b^*,\eta^*)g_{i^*}}[\mu(t_{n+1})]|\leq \frac{t_{n+1}^{-1}}2$;
\i or there exists $(\bar a,\bar b,\bar \eta,\bar i)\in \Omega'_{t_{n+1}}\times\Pg{1,\ldots,m}$ such that 
$$
|\L_{U(a^*,b^*,\eta^*)g_{i^*}}[\mu(t_{n+1})]|\leq (1-h)|\L_{U(\bar a,\bar b,\bar \eta)g_{\bar i}}[\mu(t_{n+1})]|. 
$$
\ei
\ei
}}

\medskip

The proof of convergence of $\mu(t)$ to $\mathcal{Z}$ is obtained by combining the proof of Theorem \ref{t-main} in Section \ref{s-proof} with the proof of the finite-dimensional sparse Jurdjevic--Quinn stabilization method with hysteresis given in \cite{finiteJQ}. We do not provide details.

\section{Application to crowd models}\label{sec:mas}
In this section, we give some relevant models to which Theorem~\ref{t-main} can be applied. 
Control problems for equations of the form~\eqref{e-PDEt} arise naturally when studying large crowds of interaction agents. 
Consider a system of $N$ interacting agents in which the dynamics of the state $x_i \in \R^d$ of the $i$-th agent are influenced by the state of the other $N-1$ agents, according to the time evolution
\begin{equation}\label{eq:multi}
 \dot x_i  = \frac{1}{N} \sum_{j\neq i } F(x_i,x_j), \qquad  i=1,\ldots,N,
\end{equation}
where $F\in \Lip(\R^d \times \R^d,\R^d)$ represents interaction rules, that are the same for any pair of agents.
When the number $N$ of agents is large, it is often convenient to describe the evolution of the system as a mean-field equation. 
In the mean-field limit, when $N \to +\infty$, the evolution of the mass of the agents $\mu \in \P$ is described by~\eqref{e-pdebase} with
\begin{equation}\label{eq:fmu}
f[\mu](x) = \int F(x,y) \, d\mu(y).
 \end{equation}
Indeed, to derive the mean-field model~\eqref{e-pdebase} from the finite-dimensional multi-agent models~\eqref{eq:multi}, it suffices to consider the empirical measure $\mu(t) = \frac{1}{N}\sum_{i=1}^N \delta_{x_i(t)}$.

We consider then the controlled version of the multi-agent system~\eqref{eq:multi}, given by
\begin{equation}\label{eq:multicontrol}
 \dot x_i  = \frac{1}{N} \sum_{j\neq i } F(x_i,x_j) + u_i g_i(x_1,\ldots,x_N), \qquad i=1,\ldots,N,
\end{equation}
for some Lipschitz vector field $(g_1,\ldots,g_N)$ on $(\R^d)^N$ and controls $(u_1,\ldots,u_N)$ in some subset of $(\R^d)^N$.
In the case in which the control vector field is defined, for every agent $i$, only via the interaction between the other agents and the action of the control is the same on any agent, namely if there exists an interaction kernel  $G \in \Lip(\R^d \times \R^d,\R^d)$ such that
$$
g_i(x_1,\ldots,x_N)= \frac{1}{N} \sum_{j\neq i } G(x_i,x_j) , 
$$
for every $i =1,\ldots, N$, and if $u_i=u_j$ for all $i,j \in \{1,\ldots, N\}$, then we can consider the limit of~\eqref{eq:multicontrol} as $N\to +\infty$, which gives the mean-field equation~\eqref{e-PDE} with~\eqref{eq:fmu} and
\begin{equation}\label{eq:gmu}
 g[\mu](x) = \int G(x,y) \, d\mu(y).
\end{equation}
The controllability problem is then the following: given an initial measure $\mu_0$ and a final measure $\mu_1$, find a suitable control function $(t,x) \mapsto u(t,x)$ steering the system~\eqref{e-PDE} from $\mu_0$ to $\mu_1$. 
We refer to~\cite{controlKCS} for a first result on the control of a mean-field equation of the form~\eqref{e-PDE} with constraint {\bf (U)} and {\bf (U')}.  In particular the paper focuses on the controlled version of the kinetic Cucker--Smale system introduced in~\cite{ha-tad} with constant $g$, and the existence of a control steering the system to a neighborhood of a Dirac measure is proved. 

Existence and uniqueness for the mean-field equation~\eqref{e-PDE} when the vector fields are given by~\eqref{eq:fmu} and~\eqref{eq:gmu} are ensured by Theorem~\ref{t-esistenza} provided that $F(x,y)$ and $G(x,y)$ have compact support. Indeed if $F(x,y)$ and $G(x,y)$ have compact support, then the vector fields $f[\mu]$ and $g[\mu]$ satisfy Assumptions {\bf (H)}.

Multi-agent models with a compactly supported interaction potential are sometimes called ``bounded confidence'' or homophilous models. The idea is that the agents interact only with the ones having closer states. 
This kind of interaction is used, for instance, to model opinion formation in first-order systems. 
One of the most influential models in opinion formation is, indeed, the Bounded Confidence Model by Hegselmann and Krause~\cite{HK} (see also \cite{BHT}). The main feature of this model is that the interaction is zero when the distance between two opinions is larger than a certain threshold:
$$
F(x_i,x_j) = 
\begin{cases}
(x_j-x_i) & \mbox{ if } |x_{i}-x_{j}| \leq 1,\\
0& \mbox{otherwise}.
\end{cases}
$$ 
It has been proved in \cite{BHT} that, for almost every initial configuration, the opinions converge asymptotically to clusters. 
In particular, the system does not reach global consensus in general.
Since the right-hand side is discontinuous with respect to the state variable, for some configurations, the system has no unique solution, 
hence we consider the more general first-order consensus model
\begin{equation}
\dot x_{i} = \frac{1}{N} \sum_{j\neq i} \phi(x_{j}-x_{i})(x_j-x_i)  \qquad  i=1,\ldots, N,
\label{eq:HKcont}
\end{equation}
where the function $\phi$ is defined by 
\begin{equation}
\phi(x)=\begin{cases}
1&\mbox{ if } |x|<1,\\
-\frac{|x|}{\eps}+1+\frac{1}{\eps}&\mbox{~~if~~}|x|\in[1,1+\eps],\\
0&\mbox{ if }|x |>1+\eps,
\end{cases}
\end{equation}
for some small $\eps>0$. This is a variant of the Hegselmann--Krause model, in which the Lipschitz property of $\phi$ ensures existence and uniqueness of solutions of~\eqref{eq:HKcont}. 
Therefore the associated vector field for the mean-field equation~\eqref{e-PDE} is
\begin{equation}\label{eq:fHK}
 f[\mu](x) = \int \phi(y-x)(y-x)\, d\mu(y).
\end{equation}
The kinetic version of the Hegselmann--Krause model has been first studied in~\cite{CFT12} for discrete-time dynamics. Existence of solutions has been first proved in~\cite{BCL09} for 
 a general bounded decreasing $\phi(x)$
such that $|x \phi'(x)|\leq \phi(x)$. 
Moreover, if $\phi(x)$ is everywhere nonzero, then the system converges unconditionally to consensus, meaning that for every $\mu_0$ the solution $\mu(t)$ converges asymptotically to a Dirac mass.  
If $\phi(x)$ is compactly supported, as in our case, however, then
the large time behavior of the dynamics  is not yet completely understood and a precise description of the asymptotic dynamics is, in general, a hard task. As in the finite-dimensional analogue, generically the solution $\mu(t)$ converges to a finite sum of Dirac deltas, representing the clusters of opinion, but sufficient conditions for global consensus are still unknown. 
Theorem~\ref{t-main} provides then a useful tool to establish convergence to global consensus. 

Here we consider the controlled kinetic Hegselmann--Krause model in dimension $d=1$ with drift vector field given by~\eqref{eq:fHK}, control vector field $g=1$, and Lyapunov function
$$
V[\mu] = \int x^2 \, d\mu(x).
$$
Then 
$\L_{ug} V[\mu] = 0$ for every 
$u \in \mathcal{U} \Leftrightarrow \mu = \delta_0.$ 
Since $\L_f [\delta_0]=0$, it follows that $\Z=\{\delta_0\}$. 

We claim that if $\mu  = \frac{1}{2}(\delta_x + \delta_y)$ for some $x,y \in \R$, then $\L_fV[\mu] \leq 0$.
The same statement is valid for the Lie derivative of $V$ along any combination of Dirac 
$\frac{1}{N}\sum_{i=1}^N \delta_{x_i}$ for some ${x_i} \in \R$, 
but for the sake of readability let us prove it for the sum of two Dirac masses.
Noting that 
$$
f\left[\frac{\delta_x+\delta_y}{2}\right](z) = \frac{1}{2} \phi(x-z) (x-z) +\frac{1}{2} \phi(y-z) (y-z),
$$
we have
\begin{align*}
\L_fV[\mu] &= \frac{1}{2} \frac{d}{dt}_{|_{t=0}} V(e^{tf}[\delta_x]) + \frac{1}{2} \frac{d}{dt}_{|_{t=0}} V(e^{tf}[\delta_y])\\
&= x  f\left[\frac{\delta_x+\delta_y}{2}\right](x) + y   f\left[\frac{\delta_x+\delta_y}{2}\right](y)\\
& = \frac{1}{2} \phi(y-x)(y-x)  x +\frac{1}{2} \phi(x-y)(x-y)   y\\
& = -\frac{\phi(x-y)}{2}\|x-y\|^2 \leq 0. 
\end{align*}
Then, using the continuity conditions~\eqref{e-contf} and the density of the sum of Dirac deltas in $\mathcal{P}_c(\R)$, one can extend the estimate on $\mathcal{P}_c(\R)$.

In particular this system fits into the framework of Theorem~\ref{t-main}, which thus provides the existence of a control strategy concentrating the mass at $0$, in other words, steering the system to global consensus.

Theorem~\ref{t-main} also gives an explicit construction of a control achieving consensus. Assume that at a certain switching time, say $t$, the solution is $\mu(t) = \frac{1}{2}\chi_{[-1,1]}$.  
In this case, we can write explicitly the slope function 
\begin{multline*}
s_t(a,b,\eta) = |\L_{U(a,b,\eta)g[\mu]}V[\mu]| 
=  |\frac{d}{dt}_{|_{t=0}} V(e^{t U(a,b,\eta)} \mu)| 
= 2|\int_{\R} x U(a,b,\eta)(x) d\mu(x)|\\
= 2 \left|\frac{1}{\eta} \int_{a-\eta}^a x(x-a+\eta) d\mu(x) + \int_a^b x d\mu(x) +\frac{1}{\eta} \int_b^{b+\eta} x(- x+b + \eta) d\mu(x) \right|   ,
\end{multline*}
where $b-a+2\eta \leq c$.
For $\eta$ large, the biggest contribution is given by the second integral term 
$$
\int_a^b x \, d\mu(x) = \frac{1}{4}\left(\min(b,1)^2 - \max(a,-1)^2\right).
$$
If $c<2$, then the control set $\omega$ cannot cover the whole support of $\mu$ and it will be close to $-1$ or $1$. The action of the control steers the mass in the region $\omega$ toward $0$ breaking the symmetry of the measure $\mu$. It may happen therefore that the region $\omega$ will lose the optimality of the slope function and the control will switch to another region on the opposite side. In general, if the measure $\mu$ is symmetric with respect to the origin, then the control may chatter. This is the rationale for the introduction of an hysteresis parameter $h$: the control acts on a set and holds it also sometimes after losing optimality in order to prevent high-oscillating controls.

\begin{figure}[ht]
\begin{center}
\includegraphics[width=13cm]{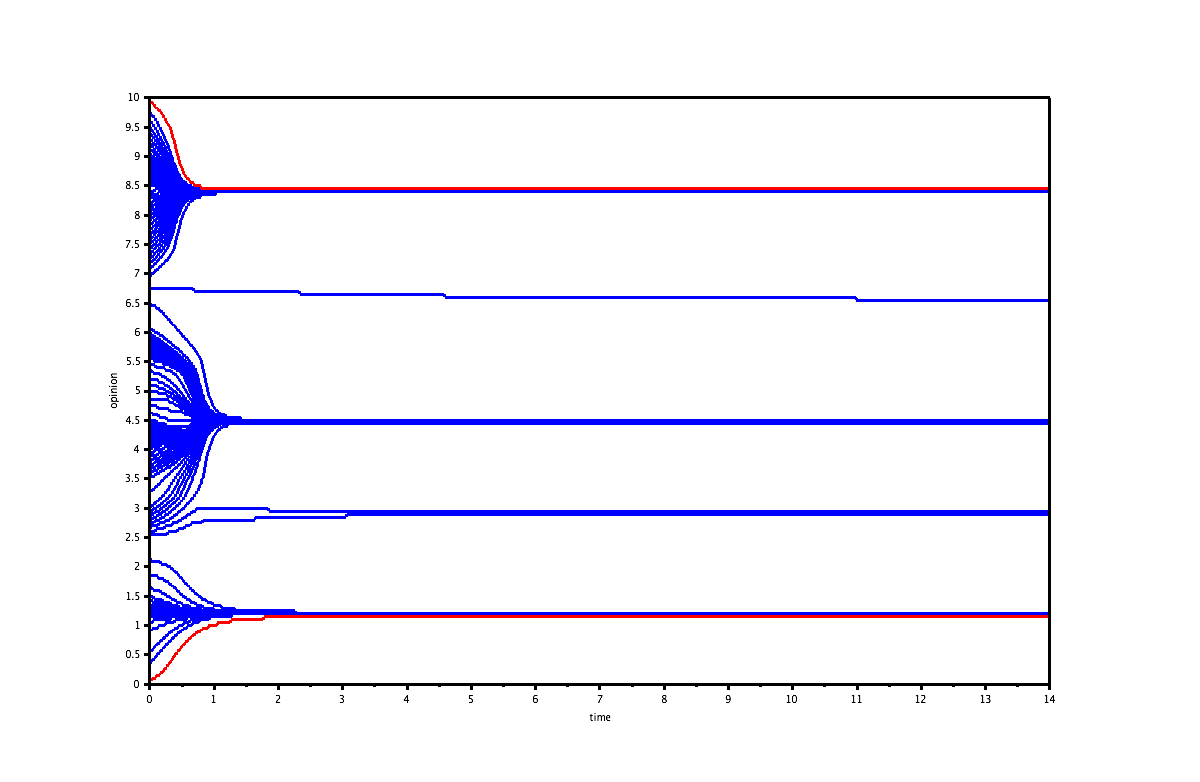}
\caption{Uncontrolled kinetic Hegselmann--Krause model. Blue lines represent the evolution $\mu(t)$. The evolution of $\max$ and $\min$ of the support of $\mu(t)$ are represented with red lines.}
\label{fig:uncontrolled}
\end{center}
\end{figure}

Here we present numerical simulations for this system. We consider an initial data $\mu_0$ randomly distributed on the interval $[0,10]$ and we apply the control given in Theorem~\ref{t-main} with three different choice of the hysteresis parameter $h$. In Figure~\ref{fig:uncontrolled} the free evolution of the system, i.e. with $u=0$. Notice that the solution tends to a finite combination of Dirac deltas representing clusters. The action of Theorem~\ref{t-main}, with the variance as Lyapunov function, is represented in Figure~\ref{fig:controlled}. In this case the whole mass tends to a single Dirac delta, representing consensus. 

\newcommand{\larghezza}{\textwidth}
\begin{figure}[ht!]
\begin{subfigure}[ht]{\larghezza}
\includegraphics[width=\larghezza]{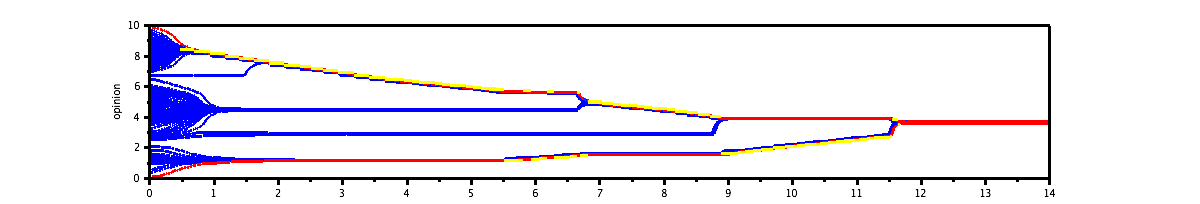}
\caption{$h=  0.9$}
\end{subfigure}

\begin{subfigure}[ht]{\larghezza}
\includegraphics[width=\larghezza]{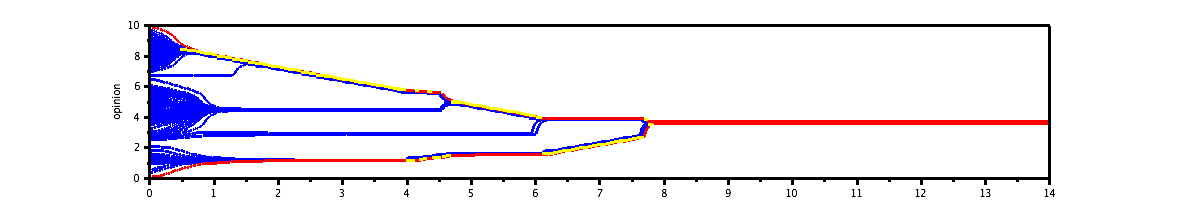}
\caption{$h =  0.5$}
\end{subfigure}

\begin{subfigure}[ht]{\larghezza}
\includegraphics[width=\larghezza]{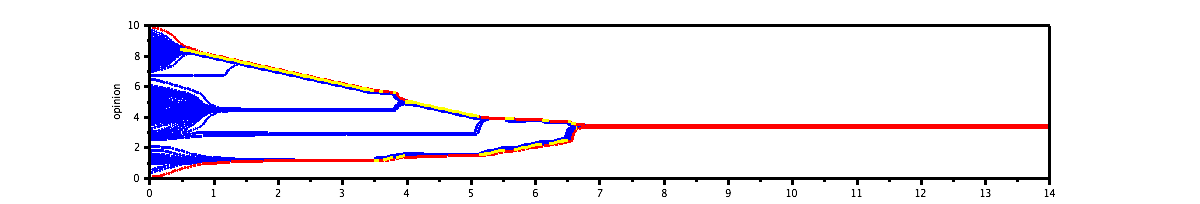}
\caption{$h =  0.2$}
\end{subfigure}

\caption{Application of the control strategy of Theorem~\ref{t-main} for different values of the hysteresis parameter $h$. The yellow region is the controlled area.}
\label{fig:controlled}
\end{figure}


Chattering may usually happen when dealing with \emph{sparse} controls designed with optimality criteria. Sparsity of the control in finite dimension is usually coded in terms of control acting on the smallest number of components/agents and the term sparsity comes from the fact that the control operator $g_i$ in~\eqref{eq:multicontrol} is a sparse vector. This notion has been introduced in~\cite{CS-control, CS-control-1} for second-order alignment systems (see also~\cite{WCB} for the controllability via leader of the Hegselmann--Krause finite-dimensional model).

In the infinite-dimensional framework, the assumption that agents are indistinguishable is crucial for defining mean-field limits; therefore, the notion of componentwise sparsity loses its sense. The infinite-dimensional analogue of componentwise sparsity is the sparsity population constraint {\bf (U')}. 
In Section~\ref{s-concentration}, we have described the issues arising from this definition.

The sparsity space constraint given in {\bf (U)} is, on the other hand, the natural definition of sparsity for mean-field equation of the form~\eqref{e-PDE}. The constraint is in some sense geographical, since the control can act only on a region of the space. 
The finite-dimensional analogue of such a control is the so-called \emph{decentralized control}. A decentralized control acts based on partial information on the agents inside a certain neighborhood of the controlled ones. The decentralized control for multi-agent systems is a well-established topic, we refer for instance to the works~\cite{dGJ,TJP07,ZJP07} for decentralized consensus algorithms, see also~\cite{BFK15} for a recent result 
with $L^\infty$ constraints on the control of the form~\eqref{e-strength}.


\section{Conclusion and open questions}
In this paper, we have generalized the classical Jurdjevic--Quinn stabilization method to infinite-dimensional control systems described by transport partial differential equations with non-local terms. Such equations arise in crowd models that are mean-field limits of particle systems for a finite number of agents: for this reason, it is natural to require some sparsity constraint to the control.

We established a mean-field Jurdjevic--Quinn stabilization method under the sparsity constraint {\bf (U)}: the control acts on a small set of the configuration space, with a bounded strength.

Improving Theorem \ref{t-main} in the original spirit of Jurdjevic and Quinn may be done in several ways: on the one side, by reducing the target goal to the largest subset of $\mathcal{Z}$ that is invariant under the uncontrolled dynamics $f[\mu]$; on the other side, by reducing the target by imposing zero higher-order derivatives, i.e., when $Z$ is defined by
\bqn
\mathcal{Z}=\Pg{\mu\in\PR\ \mid\ \L_fV[\mu]=\L_{f}^k \L_{ug_i}V[\mu]=0 \quad \forall u\in\Lip(\R^d,\R),\ i=1,\ldots,m,\ k\in\N}.
\eqnn
Dealing with iterated Lie derivatives, that is, with Lie brackets, in the kinetic setting is an open perspective.

Addressing more general systems than those presented in Section~\ref{sec:gen} is also of great interest. In particular, it would be interesting to develop similar approaches to enforce stabilization of a transport equation to a specific set $\mathcal{Z}$ of configurations, such as steady-states or periodic trajectories. In this spirit, a remarkable result for describing cell migrations is given in \cite{natalini}, where steady-states are ``rosettes'', that are symmetric configurations of cells leading to emergence of specific macroscopic structures.

\section*{Acknowledgement}
The second author acknowledges the support of the  NSF grant \#1107444 (KI-Net). The third author acknowledges the support of the ANR project CroCo ANR-16-CE33-0008. The last author acknowledges the support of the ANR project Finite4SoS ANR-15-CE23-0007-01 and of the project FA9550-14-1-0214 of the EOARD-AFOSR.


\begin{thebibliography}{99}

\bibitem{agra-capo} \auth{A.A. Agrachev, M. Caponigro}, \tit{Controllability on the group of diffeomorphisms}, jou{Annales de l'Institut Henri Poincare/Analyse non lin\'eaire} 26 (6), \pp{2503--2509}, 2009

\bibitem{ambrosio} \auth{L. Ambrosio}, \tit{Transport equation and Cauchy problem for BV vector fields}, \jou{Inventiones mathematicae}, v. 158 (2), \pp{227--260}, 2004.

\bibitem{ambrosio-gangbo} \auth{L. Ambrosio, W. Gangbo}, \tit{Hamiltonian ODEs in the Wasserstein Space of Probability Measures}, \jou{Communications on Pure and Applied Mathematics}, Volume 61, Issue 1, \pp{18--53}, 2008.

\bibitem{arg-tre} \auth{S Arguillere, E Tr\'elat}, \tit{Sub-Riemannian structures on groups of diffeomorphisms}, \jou{Journal of the Institute of Mathematics of Jussieu}, 35 pages, to appear.

\bibitem{BB-11} \auth{N. Bellomo, A. Bellouquid}, \tit{On the modeling of crowd dynamics: Looking at the beautiful shapes of swarms}, \jou{Networks and Heterogeneous Media}, 6 \pp{383--399}, 2011.

\bibitem{BHT}
\auth{V. Blondel, J. Hendrickx, and J. Tsitsiklis}, \tit{Continuous-time average-
preserving opinion dynamics with opinion-dependent communications}, \jou{SIAM Journal on Control
and Optimization}, 48(8), \pp{5214--5240}, 2010.

\bibitem{BCL09}
\auth{A. Bertozzi, J. Carrillo, and T. Laurent},  \tit{Blow-up in multidimensional aggregation equations with mildly singular interaction kernels}, \jou{Nonlinearity} 22, no. 3, \pp{683--710}, 2009.

\bibitem{BFK15}
\auth{M. Bongini, M. Fornasier, and D. Kalise}, \tit{({U}n)conditional consensus emergence under perturbed and decentralized feedback controls}, \jou{Discrete and Continuous Dynamical Systems}, Vol. 35, No. 5, \pp{4071--4094}, 2015

\bibitem{BCM}
\auth{F. Bullo, J. Cortes, S. Marti­nez}, \tit{Distributed control of robotic networks: a mathematical approach to motion coordination algorithms},
Princeton University Press, Princeton, 2009.

\bibitem{CFT12}
\auth{C. Canuto, F. Fagnani, and P. Tilli}, \tit{An Eulerian approach to the analysis of Krause's consensus models}, \jou{SIAM Journal on Control and Optimization} 50, no. 1,  \pp{243--265}, 2012.

\bibitem{CS-control} \auth{M. Caponigro, M. Fornasier, B. Piccoli, E. Tr\'elat}, \tit{Sparse stabilization and optimal control of the Cucker-Smale model}, \jou{Mathematical Control and Related Fields}, Issue 4, \pp{447--466}, 2013.

\bibitem{CS-control-1} \auth{M. Caponigro, M. Fornasier, B. Piccoli, E. Tr\'elat}, \tit{Sparse Stabilization and Control of Alignment Models}, \jou{Mathematical Models and Methods in Applied Sciences} 25, no. 3, \pp{521--564}, 2015.

\bibitem{finiteJQ} \auth{M. Caponigro, B. Piccoli, F. Rossi, E. Tr\'elat}, \tit{Sparse Jurdjevic--Quinn stabilization of dissipative systems}, preprint hal-01397843, submitted.


\bibitem{CKJRF}
\auth{I.D. Couzin, J. Krause, R. James, G.D. Ruxton, N. Franks}, \tit{Collective memory and spatial sorting in animal groups}, \jou{J Theor Biol}, 218, 2002.

\bibitem{CFP11}
\auth{E. Cristiani, P. Frasca, B. Piccoli}, \tit{Effects of anisotropic interactions on the structure of animal groups}, \jou{Journal of mathematical biology} 62, no. 4, \pp{569--588}, 2011.


\bibitem{CPT-14}
\auth{E. Cristiani, B. Piccoli, A. Tosin}, \tit{Multiscale Modeling of Pedestrian Dynamics}, Springer MS \& A: Modeling, Simulation and Applications, 2014.



\bibitem{dGJ}
\auth{M.C. De Gennaro and A. Jadbabaie}, \tit{Decentralized control of connectivity for multi-agent systems}, In \jou{Proceedings of the 45th IEEE Conference on Decision and Control}, \pp{3628--3633}, 2006.

\bibitem{natalini} \auth{E. Di Costanzo, R. Natalini, L. Preziosi}, \tit{A hybrid mathematical model for self-organizing cell migration in the zebrafish lateral line}, \jou{J. Math. Biol.} 71, \pp{171--214}, 2015.

\bibitem{FPR-14}
\auth{M. Fornasier, B. Piccoli, F. Rossi}, \tit{Mean-field sparse optimal control}, \jou{Philosophical Transaction Royal Society} A, 372, 2014.

\bibitem{G-08}
\auth{I. Giardina},  \tit{Collective behavior in animal groups: theoretical models and empirical studies}, \jou{Human Frontier Science Program Journal}, (205–219), 2008.

\bibitem{paola} \auth{P. Goatin, F. Rossi},  \tit{A traffic flow model with non-smooth metric interaction: well-posedness and micro-macro limit}, \jou{Comm. Math. Sciences}, to appear.



\bibitem{ha-tad} \auth{S.-Y. Ha, E. Tadmor}, \tit{From particle to kinetic and hydrodynamic description of flocking}, \jou{Kinetic and Related Methods}, v. 1 (3), \pp{415--435}, 2008.

\bibitem{HK}
\auth{R. Hegselmann, U. Krause}, \tit{Opinion dynamics and bounded confidence models, analysis, and simulation}, \jou{Journal of Artificial Societies and Social Simulation}, 5(3), 2002.

\bibitem{JLM-03}
\auth{Ali Jadbabaie, Jie Lin, A Stephen Morse}, \tit{Coordination of groups of mobile autonomous agents using nearest neighbor rules},  \jou{Automatic Control, IEEE Transactions on}, 48(6) \pp{988--1001}, 2003.

\bibitem{HG-02}
\auth{R. Hegselmann, U. Krause}, \tit{Opinion dynamics and bounded confidence: models, analysis and simulation}, \jou{Journal of Artificial Societies and Social Simulation}, 5(3), 2002.

\bibitem{HFV-00}
\auth{D. Helbing, I. Farkas, T. Viscek}, \tit{Simulating dynamical features of escape panic}, \jou{Nature}, 407, \pp{487--490}, 2000.

\bibitem{HPS15}
\auth{M. Herty, L. Pareschi, S. Steffensen}, \tit{Mean-field control and Riccati equations}, \jou{Networks
and Heterogeneous Media}, 10(3), \pp{699--715}, 2015.


\bibitem{JQ} \auth{V. Jurdjevic, J.P. Quinn}, \tit{Controllability and stability}, \jou{Journal of Differential Equations}, Volume 28, Issue 3, \pp{381--389}, 1978.

\bibitem{OSFM}
\auth{Reza Olfati-Saber, J Alex Fax, Richard M Murray},  \tit{Consensus and cooperation in networked multi-agent systems}, Proceedings of the IEEE, 95(1): \pp{215--233}, 2007.

\bibitem{genwass} \auth{B. Piccoli, F. Rossi}, \tit{Generalized Wasserstein distance and its application to transport equations with source}, \jou{Archive for Rational Mechanics and Analysis}, Volume 211, Issue 1, pp. 335-358, 2014.

\bibitem{pedestrian} \auth{B. Piccoli, F. Rossi}, \tit{Transport equation with nonlocal velocity in Wasserstein spaces: convergence of numerical schemes}, \jou{Acta Applicandae Mathematicae}, 124, \pp{73--105}, 2013.

\bibitem{controlKCS} \auth{B. Piccoli, F. Rossi, E. Tr\'elat}, \tit{Control to flocking of the kinetic Cucker-Smale model}, \jou{SIAM J. Mathematical Analysis} 47, no. 6, \pp{4685--4719}, 2015.


\bibitem{RME}
\auth{A. Rahmani, M. Ji, M. Mesbahi, M. Egerstedt}, \tit{Controllability of multi-agent systems from a graph-theoretic perpective}, \jou{SIAM Journal on Control and Optimization}, 48(1) \pp{162--186}, 2009.

\bibitem{TJP07}
\auth{H. Tanner, A. Jadbabaie, and G. Pappas}, \tit{Flocking in fixed and switching networks}, \jou{IEEE Transactions on Automatic Control}, 52(5), \pp{863--868}, 2007.


\bibitem{villani} \auth{C. Villani}, \tit{Topics in Optimal Transportation}, Graduate Studies in Mathematics, Vol. 58, 2003.


\bibitem{WCB}
\auth{S. Wongkaew, M. Caponigro, and A. Borzi}, \tit{On the control through leadership of the Hegselmann--Krause opinion formation model}, \jou{Mathematical Models and Methods in Applied Sciences}, Volume 25, Issue 03, \pp{565--585}, (2015). 

\bibitem{ZJP07}
\auth{M. Zavlanos, A. Jadbabaie, and G. Pappas}, \tit{Flocking while preserving network connectivity}, In \jou{Proceedings of the 46th IEEE Conference on Decision and Control}, \pp{2919--2924}, 2007.


\end{thebibliography}
\end{document}